\documentclass[12pt]{amsart}%
\usepackage{indentfirst, latexsym, bm, amsmath, amsthm, amscdx, amssymb,
mathabx,mathrsfs}
\usepackage[colorlinks, linkcolor=blue,citecolor=red]{hyperref}
\usepackage{upgreek}
\usepackage{bm}
\usepackage{url}
\usepackage{lineno}
\usepackage{amsmath}
\usepackage{amsfonts}
\usepackage{amssymb}
\usepackage{graphicx}%
\providecommand{\U}[1]{\protect\rule{.1in}{.1in}}
\theoremstyle{definition}
\theoremstyle{remark}
\numberwithin{equation}{section}
\newtheorem{theorem}{Theorem}[section]
\newtheorem*{theorem*}{Theorem}

\newtheorem{lemma}{Lemma}[section]
\newtheorem*{lemma*}{Lemma}
\newtheorem{corollary}{Corollary}[section]
\newtheorem{remark}{Remark}[section]
\newtheorem{proposition}{Proposition}[section]
\newtheorem{definition}{Definition}[section]

\begin{document}
\title[Legendrian Self-Shrinker]{Rigidity and Classification of Legendrian Self-Shrinkers}
\author[Chang]{Shu-Cheng Chang$^{1\ast}$}
\address{$^{1}$Department of Mathematics, National Taiwan University, Taipei 10617,
Taiwan and Shanghai Institute of Mathematics and Interdisciplinary Sciences,
Shanghai, 200433, China }
\email{scchang@math.ntu.edu.tw}
\author{Chin-Tung Wu$^{2\ast\ast}$}
\address{$^{2}$Department of Applied Mathematics, National Pingtung University,
Pingtung 90003, Taiwan}
\email{ctwu@mail.nptu.edu.tw }
\author{Liuyang Zhang$^{3\ast\ast\ast}$}
\address{$^{3}$Mathematical Science Research Center, Chongqing University of
Technology, 400054, Chongqing, P.R. China}
\email{zhangliuyang@cqut.edu.cn}
\author{Qiuxia Zhang$^{4}$}
\address{$^{4}$Mathematical Science Research Center, Chongqing University of
Technology, 400054, Chongqing, P.R. China}
\email{2430225058@qq.com}
\thanks{$^{\ast}$Research supported in part by Startup Foundation for Advanced Talents
of the Shanghai Institute for Mathematics and Interdisciplinary Sciences
(No.2302-SRFP-2024-0049). $^{\ast\ast}$Research supported in part by NSTC
grant 113-2115-M-153-003, Taiwan. $^{\ast\ast\ast}$ Research supported in part
by the Scientific and Technological Research Program of Chongqing Municipal
Education Commission (No. KJQN202201138); Startup Foundation for Advanced
Talents of Chongqing University of Technology (No. 2022ZDZ019)}
\maketitle

\begin{abstract}
In this article, we first classify Legendrian self-shrinkers in $\mathbb{R}%
^{3}$ and $\mathbb{R}^{5}$. We then prove a Legendrian rigidity theorem,
which can be regarded as an analogue of the result of Li-Wang \cite{lw}. More
precisely, let $F(\Sigma)\subset\mathbb{R}^{5}$ be an orientable Legendrian
self-shrinker, if $\Vert A\Vert_{g}^{2}\leq2$ and the associated Legendrian
immersion $\bar{F}\subset\mathbb{R}^{4}\times\mathbb{S}^{1}$ is compact, then
$\bar{F}$ must be a flat minimal generalized Legendrian Clifford torus in
$\mathbb{S}^{5}$, whose cone $\mathcal{C}(\bar{F}(\Sigma))$ is the
Harvey-Lawson special Lagrangian cone in $\mathbb{C}^{3}$.

\end{abstract}

\section{Introduction}

It was proven by Haskins in \cite{h} that $L$ is Legendrian in the standard
Sasakian sphere $\mathbb{S}^{2n+1}$ if and only if the cone $\mathcal{C}(L)\setminus
\{0\}$ is Lagrangian in $\mathbb{C}^{n+1}$. In fact, there exists a one-to-one
correspondence between minimal Lagrangian cones in $\mathbb{C}^{n+1}$ and
minimal Legendrian submanifolds in $\mathbb{S}^{2n+1}$. Furthermore,
$\mathbb{S}^{2n+1}$ can be viewed as the compact contactization of the
standard K\"{a}hler projective space $\mathbb{CP}^{n}$. In particular, the
projection of the mean curvature vector field of a Legendrian submanifold $L$
in the contactization of a K\"{a}hler manifold $Z$ coincides with the mean
curvature of the projected Lagrangian submanifold in $Z$.

On the other hand, any closed $n$-manifold $L$ can be embedded in the contact
Euclidean space $\mathbb{R}^{2n+1}$, and any two such embeddings are isotopic
for $n\geq2$. However, the Legendrian isotopy problem remains open due to the
involvement of contact invariants. For $n=1$, Eliashberg-Fraser \cite{ef} and
Colin-Giroux-Honda \cite{cgh} demonstrated that there are only finitely many
Legendrian knot types, which are classified by the Thurston-Bennequin
invariant, a classical rotation number, and the knot type. However, for $n>1$,
Ekholm-Etnyre-Sullivan \cite{ees} showed that there exists an infinite family
of Legendrian embeddings of the $n$-sphere (or $n$-torus) into $\mathbb{R}%
^{2n+1}$ that are not Legendrian isotopic, even though they share the same
rotation number and Thurston-Bennequin invariant.

To address this difficulty, Smoczyk \cite{s1} first introduced the Legendrian
mean curvature flow, which preserves the Legendrian condition if the ambient
manifold is an $\eta$-Einstein Sasakian manifold $(M^{2n+1},\Phi,\xi,\eta,g)$:%
\begin{equation}%
\begin{array}
[c]{c}%
\frac{d}{dt}F_{t}=H-2\alpha\xi
\end{array}
\label{c}%
\end{equation}
where $H$ is the Legendrian mean curvature vector and $\alpha$ is the Legendre
angle of $L_{t}=F_{t}(L)$, normalized as in \cite{s1} by
\begin{equation}%
\begin{array}
[c]{c}%
F_{t}^{\ast}\Omega^{T}=e^{-i\alpha}Vol_{L_{t}},\text{ \ \ equivalently \ }%
H=\Phi\nabla\alpha,
\end{array}
\label{angle}%
\end{equation}
where $\Omega^{T}$ is the transverse holomorphic volume form of the ambient
Sasakian manifold and $Vol_{L_{t}}$ the induced volume form (see Section 2).

Both the sign and the factor $2$ in (\ref{c}) are dictated by the Sasakian
normalization $d\eta=2g(\cdot,\Phi\cdot)$ used throughout this paper, see
(\ref{eq2.3a}). Indeed, along a Legendrian immersion the normal bundle splits
as $NL=\Phi TL\oplus\mathbb{R}\xi$, so that every normal variation field can be
written as $V=\Phi W+h\xi$ with $W\in TL$ and a smooth function $h$, and
Cartan's formula together with (\ref{eq2.3a}) gives
\[%
\begin{array}
[c]{c}%
\frac{\partial}{\partial t}\left(  F_{t}^{\ast}\eta\right)  =dh+2W^{\flat}.
\end{array}
\]
Hence the Legendrian condition $F_{t}^{\ast}\eta=0$ is preserved if and only if
$W=-\frac{1}{2}\nabla h$; since $H=\Phi\nabla\alpha$, this forces
$h=-2\alpha$ up to an additive constant, which is precisely (\ref{c}). This is
Lemma 3.1 of \cite{s1} in our notation, and (\ref{c}) is literally Smoczyk's
Legendrian mean curvature flow $\frac{d}{dt}F=\vec{H}-2\alpha X_{\lambda}$.
In particular, replacing $2\alpha\xi$ by $\alpha\xi$ in (\ref{c}) would destroy
the Legendrian condition unless $\alpha$ is constant, i.e. unless $L_{t}$ is
minimal Legendrian.

Smoczyk proved that any Legendrian
unknot with the same rotation number and Thurston-Bennequin invariant is
Legendrian isotopic to a closed Legendrian geodesic in a tight Sasakian
3-manifold with nonnegative $\eta$-Einstein curvature. Recently, the first two
authors \cite{chw} proved the existence of long-time solutions and asymptotic
convergence along the Legendrian mean curvature flow in higher-dimensional
$\eta$-Einstein Sasakian $(2n+1)$-manifolds, under suitable stability
conditions stemming from the Thomas-Yau conjecture \cite{TY01}.

However, singularities are inevitable in the Legendrian mean curvature flow,
and understanding how these singularities form is crucial for studying the
isotopy problem through the flow method. Addressing these issues requires a
deep understanding of the Legendrian mean curvature flow and the types of
singularities that may arise in finite time. In fact, singularities in
Legendrian mean curvature flow are generally modeled locally by self-similar
soliton solutions, such as Legendrians in the standard contact Euclidean space
$(\mathbb{R}^{2n+1}, \Phi, \xi, \eta, g)$, which evolve via rescaling or
translation under the mean curvature flow.

In the work of Chang-Wu-Zhang \cite{cwz}, the authors investigated Type-I
singularities of the Legendrian mean curvature flow through blow-up analysis.
If the singularity is Type-I, then there exists a subsequence such that the
space-time track of the smooth flow, under parabolic dilation, exists smoothly
on $(-\infty,0)$ and satisfies%
\[%
\begin{array}
[c]{c}%
H-2\alpha\xi=\frac{1}{2s}X^{\bot}%
\end{array}
\]
on $\mathcal{F}_{\infty}$ for $-\infty<s<0$, where
\[%
\begin{array}
[c]{c}%
X:=\sum_{i=1}^{n}\left(  x_{i}\frac{\partial}{\partial x_{i}}+y_{i}%
\frac{\partial}{\partial y_{i}}\right)  +2z\frac{\partial}{\partial z}%
\end{array}
\]
is the infinitesimal generator of the contact dilations $\delta_{\lambda
}(x,y,z)=(\lambda x,\lambda y,\lambda^{2}z)$ of $\mathbb{R}^{2n+1}$, see
Definition \ref{dilation}. That is, it is the self-shrinker solution
\[
F(x,s)=\delta_{\sqrt{-s}}\left(  F(x,-1)\right)
\]
of the Legendrian mean curvature flow.

With a slight change of notations from \cite{cwz}, we choose the standard
contact form and associated metric in $\mathbb{R}^{2n+1}$ (see Section 2 for
more details). In this context, we define the following:

\begin{definition}
Let $F:L^{n}\rightarrow(\mathbb{R}^{2n+1},\Phi,\xi,\eta,g)$ be an
$n$-dimensional Legendrian submanifold in the contact Euclidean space
$(\mathbb{R}^{2n+1},\Phi,\xi,\eta,g)$. We say that the Legendrian immersed
manifold $L^{n}$ is a self-shrinker or a self-expander if it satisfies the
following quasilinear elliptic system:
\begin{equation}%
\begin{array}
[c]{c}%
H-2\alpha\xi\mathbf{=}aX^{\perp},
\end{array}
\label{CCC}%
\end{equation}
where $a<0$ or $a>0$, respectively. Here, $H$ denotes the Legendrian mean
curvature vector, $\alpha$ the Legendre angle (\ref{angle}), $X$ the generator
of the contact dilations of Definition \ref{dilation}, and $X^{\perp}$ its
projection onto the normal bundle
\[
NL:=\Phi TL\oplus\mathbb{R\xi}%
\]
of $L^{n}$.
\end{definition}

\begin{remark}
\label{rem1.2}Decomposing (\ref{CCC}) along $NL=\Phi TL\oplus\mathbb{R}\xi$
shows that it is equivalent to the pair of equations
\[%
\begin{array}
[c]{c}%
H=a\left(  X^{\perp}\right)  _{\Phi TL},\text{ \ \ }-2\alpha=az.
\end{array}
\]
We stress that the self-similarity of these solutions is with respect to the
contact dilations $\delta_{\lambda}$ of Definition \ref{dilation}, whose
generator is $X$ and not the Euclidean position vector: the Legendrian
Clifford torus of Theorem \ref{rigidity}, for instance, evolves under
(\ref{c}) as $F_{\tau}(t,s)=\left(  re^{it},re^{is},-\frac{r^{2}}%
{2}(t+s)\right)  $ with $r^{2}=r^{2}(\tau)=-8\tau$, so that the contact
direction $z$ scales like $r^{2}$, and not like $r$, as $\tau\nearrow0$.
\end{remark}

For the Lagrangian mean curvature flow, Castro-Lerma \cite{cl} and
Joyce-Lee-Tsui \cite{jlt} constructed several examples of Lagrangian
self-similar solutions and translating solitons.

In this paper, we focus on the rigidity theorem and the classification of
Legendrian self-shrinkers in the contact Euclidean five-space $\mathbb{R}^{5}%
$. Specifically, in Section 3, we construct numerous examples of self-similar
solutions to the Legendrian mean curvature flow:

\begin{theorem}
\label{Ex} Write $w_{j}=r_{j}(s)e^{i\phi_{j}(s)},r_{j}=|w_{j}|$ and $\phi
=\phi_{1}+\phi_{2}$, for functions $r_{j}:I\rightarrow(0,\infty)$ and
$\phi_{1},$ $\phi_{2},$ $\phi:I\rightarrow\mathbb{R}$ or $\mathbb{R}{/2}%
\pi\mathbb{Z}$. Fix $s_{0}\in{I}$. Define $u:I\rightarrow\mathbb{R}$ by
\begin{equation}%
\begin{array}
[c]{c}%
u(s)=2\int_{s_{0}}^{s}r_{1}(t)r_{2}(t)\cos(\alpha(t)+\phi(t))dt.
\end{array}
\label{eq15a}%
\end{equation}
For the above $u$, we have $r_{j}^{2}(s)=\rho_{j}+\lambda_{j}u(s),$
$j=1,2,s\in I,$ $\rho_{j}=r_{j}^{2}\left(  s_{0}\right)  $. Define a degree
$2$ polynomial $Q(u)$ by $Q(u)=\prod_{j=1}^{2}\left(  \rho_{j}+\lambda
_{j}u\right)  $. Suppose that%
\begin{equation}
\left\{
\begin{array}
[c]{l}%
\frac{dw_{1}}{ds}=\lambda_{1}e^{-i\alpha(s)}\overline{w}_{2},\text{
\ }\frac{dw_{2}}{ds}=\lambda_{2}e^{-i\alpha(s)}\overline{w}_{1},\\
\frac{d\phi_{j}}{ds}=-\frac{\lambda_{j}}{\rho_{j}+\lambda_{j}u}%
Q(u)^{1/2}\sin(\alpha+\phi),\text{ }j=1,2,\\
\frac{dz}{ds}=\frac{C}{2}Q(u)^{1/2}\sin(\alpha+\phi),\end{array}
\right.  \label{eq16a}%
\end{equation}
holds in $I$. Then the submanifold $\tilde{L}$ in $\mathbb{C}{^{2}}{{\times
}\mathbb{R}}$ given by
\begin{equation}%
\begin{array}
[c]{c}%
\tilde{L}=\left\{  \left(  x_{1}w_{1}(s),x_{2}w_{2}(s),z(s)\right)  :s\in
I,x_{j}\in\mathbb{R},\sum_{j=1}^{2}\lambda_{j}x_{j}^{2}=C\right\}  ,
\end{array}
\label{eq17a}%
\end{equation}
is Legendrian. Furthermore, if it holds
\begin{equation}%
\begin{array}
[c]{c}%
\frac{d\alpha}{ds}=-\frac{a}{4}Cr_{1}r_{2}\sin(\alpha+\phi)
\end{array}
\label{eq19a}%
\end{equation}
in $I$, then $\tilde{L}$ satisfies \eqref{CCC} and is a self-similar solution
to the Legendrian mean curvature flow.
\end{theorem}

In the one-dimensional case the classification is complete; in Section 4 we
construct the Legendre analogue of the Abresch-Langer curve and prove

\begin{theorem}
\label{A-L} Let $\tilde{\gamma}(t):I\rightarrow(\mathbb{R}^{3},\Phi,\xi
,\eta,g)$ be an immersed curve as
\[%
\begin{array}
[c]{c}%
\tilde{\gamma}(t)=\left(  y(t)\cos\alpha(t)+x(t)\sin\alpha(t),x(t)\cos
\alpha(t)-y(t)\sin\alpha(t),\frac{B}{2}\alpha(t)\right)
\end{array}
\]
with smooth functions $x(t),$ $y(t),$ $\alpha(t)$ solved by the ODE system:
\[%
\begin{cases}
\frac{d}{dt}x(t)=x(t)y(t)\\
\frac{d}{dt}y(t)=B-x^{2}(t)\\
\frac{d}{dt}\alpha(t)=x(t),
\end{cases}
\]
where $B$ is a positive constant. Then $\tilde{\gamma}(t)$ is a Legendrian
curve with Legendre angle $\alpha$ and curvature $\kappa=\frac{2}{B}x$, and it
solves
\[
\kappa\hat{N}-2\alpha(t)\xi=-\frac{4}{B}X^{\perp}.
\]
In particular, when $B=1$, $\tilde{\gamma}$ is the self-shrinker of the
Legendre curve shortening flow
\[
\tilde{\gamma}_{t}=\kappa\hat{N}-2\alpha(t)\xi,
\]
where $\hat{N}:=\Phi\frac{\dot{\tilde{\gamma}}(t)}{\Vert\dot{\tilde{\gamma}%
}(t)\Vert_{g}}$ is the normal vector field on $\tilde{\gamma}$.
\end{theorem}

In the two-dimensional case the partial classification result is

\begin{theorem}
\label{c-l} Let $F:M^{2}\rightarrow(\mathbb{R}^{5},\Phi,\xi,\eta,g)$ be an
immersed self-shrinker for Legendrian mean curvature flow with Legendre
angle $\alpha$ a harmonic function with respect to $g^{T}$, then $F$ is an
embedding and locally congruent to one of the following models:

(1) $\tilde{C}:\mathbb{R}^{2}\rightarrow\mathbb{C}^{2}\times\mathbb{R}$
\[%
\begin{array}
[c]{c}%
\tilde{C}(t,s)=\left(  \frac{1}{\sqrt{-2a}}e^{it},s,\frac{1}{4a}t+C_{1}\right)
,
\end{array}
\]

(2) $\tilde{T}:\mathbb{R}^{2}\rightarrow\mathbb{C}^{2}\times\mathbb{R}$
\[%
\begin{array}
[c]{c}%
\tilde{T}(t,s)=\left(  \frac{1}{\sqrt{-2a}}e^{it},\frac{1}{\sqrt{-2a}}%
e^{is},\frac{1}{4a}(t+s)+C_{2}\right)  ,
\end{array}
\]

(3) $\tilde{\Upsilon}_{\gamma}:\mathbb{R}^{2}\rightarrow\mathbb{C}^{2}%
\times\mathbb{R}$ with parameter $0<\gamma<\frac{\pi}{2}$,
\[%
\begin{array}
[c]{c}%
\tilde{\Upsilon}_{\gamma}(s,t)=\left(  \frac{-i\sin\gamma}{\sqrt{-2a}}\cosh
te^{\frac{is}{\cos\gamma}},\frac{\tan\gamma}{\sqrt{-2a}}\sinh te^{-i\cos\gamma
s},\frac{\sin\gamma\tan\gamma}{4a}s+C_{3}\right)  ;
\end{array}
\]

(4) $\tilde{\Psi}_{\nu}:\mathbb{S}^{1}\times\mathbb{R}\rightarrow
\mathbb{C}^{2}\times\mathbb{R}$ with parameter $\nu>0$,
\[%
\begin{array}
[c]{c}%
\tilde{\Psi}_{\nu}(e^{is},t)=\left(  \frac{\cosh\nu}{\sqrt{-2a}}\cos
se^{\frac{it}{\sinh\nu}},\frac{\coth\nu}{\sqrt{-2a}}\sin se^{i\sinh\nu
t},\frac{\cosh\nu\coth\nu}{4a}t+C_{4}\right)  ,
\end{array}
\]
where $C_{i},$ $i=1,2,3,4$ are constants and the first two components in
(1)--(4) are the Hamiltonian stationary Lagrangian models of Castro-Lerma
\cite{cl}. Moreover, they all satisfy $H-2\alpha\xi=8aX^{\perp}$ with $a<0$.
\end{theorem}

Applying the Lagrangian projection (see Proposition \ref{P51}) and Theorem
\ref{c-l}, we obtain the Legendrian analogue of Li-Wang's rigidity result in
\cite{lw}:

\begin{theorem}
\label{rigidity} Let $F:\Sigma^{2}\rightarrow(\mathbb{R}^{5},\Phi,\xi,\eta,g)$
be a Legendrian self-shrinker
\[%
\begin{array}
[c]{c}%
H-2\alpha\xi=-X^{\perp},
\end{array}
\]
denote $\bar{F}$ the corresponding Legendrian immersion in $\mathbb{R}%
^{4}\times\mathbb{S}^{1}$ and $A$ the second fundamental form of $F$.

If $\bar{F}(\Sigma)$ is compact and $\Vert A\Vert_{g}^{2}\leq2$, then $\Vert
A\Vert_{g}^{2}=2$,
\[
F(t,s)=(2e^{it},2e^{is},-2(t+s)),
\]
and%
\[
\bar{F}(t,s)=(\sqrt{2}e^{it},\sqrt{2}e^{is},e^{-i(t+s)}).
\]
Moreover, after normalizing the moduli of the three complex coordinates,
$\bar{F}(\Sigma)$ is the Harvey-Lawson torus
\[
\mathbb{T}_{HL}=\left\{  \tfrac{1}{\sqrt{3}}(e^{i\vartheta_{1}},e^{i\vartheta
_{2}},e^{i\vartheta_{3}})\in\mathbb{S}^{5}:\vartheta_{1}+\vartheta_{2}%
+\vartheta_{3}\equiv0\ (\mathrm{mod}\ 2\pi)\right\}  ,
\]
the flat minimal Legendrian torus of $\mathbb{S}^{5}$ whose cone
\[
\mathcal{C}(\mathbb{T}_{HL})=\left\{  (z_{1},z_{2},z_{3})\in\mathbb{C}%
^{3}:|z_{1}|=|z_{2}|=|z_{3}|,\text{ }z_{1}z_{2}z_{3}\in\lbrack0,\infty
)\right\}
\]
is the Harvey-Lawson special Lagrangian cone \cite{hl}.
\end{theorem}

The article is organized as follows. In Section $2$, we provide a detailed
proof of the propositions concerning the adapted frames selected for the
calculations throughout the article. Additionally, we briefly introduce the
Lagrangian projection and the Legendrian lifting theorem. Section $3$ is
dedicated to the construction of various examples of self-similar solutions to
the Legendrian mean curvature flow, along with a comprehensive proof of
Theorem \ref{Ex}. In Section $4$, we focus on the proof of the classification
results in Theorems \ref{A-L} and \ref{c-l}. Finally, in the concluding part
of this section, we establish the rigidity result presented in Theorem
\ref{rigidity}.

\section{Preliminary}

This section is devoted to settings of the ambient contact space and
definition of the self-similar solution to the Legendrian mean curvature flow.

Consider the space $(\mathbb{R}^{5},\Phi,\xi,\eta,g)$ equipped with $4$
structures, which are the contact $1$-form
\[%
\begin{array}
[c]{c}%
\eta=\frac{1}{2}d{z}-\frac{1}{4}\sum_{i=1}^{2}({y_{i}}d{x_{i}}-{x_{i}}d{y_{i}%
}),
\end{array}
\]
the associated metric
\[%
\begin{array}
[c]{c}%
g=\eta\otimes\eta+\frac{1}{4}\sum_{i=1}^{2}((d{x_{i}})^{2}+(d{y_{i}})^{2})
\end{array}
\]
with $z_{j}=x_{j}+\sqrt{-1}y_{j}$, $\omega_{\mathrm{can}}=\frac{\sqrt{-1}}{8}%
(dz_{1}\wedge d\overline{z}_{1}+\cdots+dz_{n}\wedge d\overline{z}_{n})=\frac
{1}{4}\sum_{i=1}^{n}dx_{i}\wedge dy_{i}$ (normalized so that $\omega
_{\mathrm{can}}=g^{T}(\cdot,\Phi\cdot)$), the
Reeb vector field $\xi=2\frac{\partial}{\partial z}$ and the $(1,1)$-tensor
$\Phi$ expressed by the matrix
\begin{equation}%
\begin{array}
[c]{l}%
\Phi\sim\left(
\begin{array}
[c]{ccccc}%
0 & 0 & 1 & 0 & 0\\
0 & 0 & 0 & 1 & 0\\
-1 & 0 & 0 & 0 & 0\\
0 & -1 & 0 & 0 & 0\\
\frac{x_{1}}{2} & \frac{x_{2}}{2} & \frac{y_{1}}{2} & \frac{y_{2}}{2} & 0
\end{array}
\right)  .
\end{array}
\label{eq2.1}%
\end{equation}

It is evident that $(\mathbb{R}^{5},g)$ is a Riemannian manifold. In fact, we
will prove that it is a Sasakian manifold, which is equivalent to the
following condition:%
\begin{align}
\Phi^{2}+I  &  =\xi\otimes\eta,\label{eq2.2a}\\
\frac{1}{2}d\eta(\cdot,\cdot)  &  =g(\cdot,\Phi\cdot),\label{eq2.3a}\\
g(\cdot,\cdot)  &  =g(\Phi\cdot,\Phi\cdot)+\eta\otimes\eta,\label{eq2.4a}\\
\mathcal{N}_{\Phi}+d\eta\otimes\xi &  =0. \label{eq2.5a}%
\end{align}

\begin{definition}
A Sasakian manifold $(M,\Phi,\xi,\eta, g)$ is called $\eta$-Einstein if there
is a constant $K$ such that the Ricci curvature
\[
Ric=Kg+(2n-K)\eta\otimes\eta
\]
and then
\[
Ric^{T}=(K+2)g^{T}%
\]
which is called the transverse K\"{a}hler-Einstein. It is called
Sasaki-Einstein if $K=2n$. Then
\[
Ric=2ng.
\]
In case of the contact Euclidean space $\mathbb{R}^{2n+1}$, then $K=-2$ and
\[
Ric=-2g+(2n+2)\eta\otimes\eta.
\]

\end{definition}

\begin{remark}
The ambient space $(\mathbb{R}^{5},\Phi,\xi,\eta,g)$ is $\eta$-Einstein, a
property that is essential for defining the Legendrian mean curvature flow, as
proved in \cite{s1}. In the following, we denote $g^{T}=\frac{1}{4}\sum
_{i=1}^{2}((d{x_{i}})^{2}+(d{y_{i}})^{2})$ as the transverse part of $g$.
\end{remark}

Besides the contact structure itself, we shall need the scaling which is
adapted to it.

\begin{definition}
\label{dilation}For $\lambda>0$ let
\[%
\begin{array}
[c]{c}%
\delta_{\lambda}:\mathbb{R}^{2n+1}\rightarrow\mathbb{R}^{2n+1},\text{
\ }\delta_{\lambda}(x_{1},y_{1},\cdots,x_{n},y_{n},z)=(\lambda x_{1},\lambda
y_{1},\cdots,\lambda x_{n},\lambda y_{n},\lambda^{2}z)
\end{array}
\]
be the contact dilations of $(\mathbb{R}^{2n+1},\Phi,\xi,\eta,g)$, which
satisfy
\[
\delta_{\lambda}^{\ast}\eta=\lambda^{2}\eta,\text{ \ \ }\delta_{\lambda}^{\ast
}g^{T}=\lambda^{2}g^{T}.
\]
Their infinitesimal generator is the vector field
\begin{equation}%
\begin{array}
[c]{c}%
X:=\sum_{i=1}^{n}\left(  x_{i}\frac{\partial}{\partial x_{i}}+y_{i}%
\frac{\partial}{\partial y_{i}}\right)  +2z\frac{\partial}{\partial z}%
=F+z\frac{\partial}{\partial z}=F+\frac{z}{2}\xi,
\end{array}
\label{eqX}%
\end{equation}
where $F$ denotes the Euclidean position vector field. In particular
$\eta(X)=z=2\eta(F)$, whereas $X-F$ is proportional to $\xi$, so that $X$ and
$F$ have the same component in the contact distribution $\ker\eta$.
\end{definition}

\begin{remark}
It is $\delta_{\lambda}$, and not the Euclidean dilation, that is compatible
with the contact structure, and it is $\delta_{\lambda}$ that appears in the
parabolic rescaling of the Legendrian mean curvature flow (\ref{c}). This is
why the self-similar solutions (\ref{CCC}) are modelled on $X^{\perp}$ and not
on the Euclidean position field $F$, whose flow lines are not compatible with
the contact structure.
\end{remark}

A Sasakian manifold is endowed with a unique curvature tensor, which serves a
role similar to the holomorphic sectional curvature in K\"{a}hler geometry.

\begin{definition}
For any vector $X$ orthogonal to the Reeb vector field $\xi$, the sectional
curvature $K(X,\Phi X)$ denoted by $\mathcal{H}(X)$ is called $\Phi$-sectional curvature.
\end{definition}

\begin{proposition}
[Theorem 7.18 of \cite{bl}]\label{curvature} The $\Phi$-sectional curvature of
a Sasakian manifold determines the Riemannian curvature tensor completely.
\end{proposition}

\begin{remark}
\label{rem2.1} In other literature, $\eta$ form is chosen to be
\[%
\begin{array}
[c]{c}%
\tilde{\eta}=\frac{1}{2}d{z}-\frac{1}{2}\sum_{i=1}^{2}{y_{i}}d{x_{i}},
\end{array}
\]
with associated metric
\[%
\begin{array}
[c]{c}%
\tilde{g}=\tilde{\eta}\otimes\tilde{\eta}+\frac{1}{4}\sum_{i=1}^{2}((d{x_{i}%
})^{2}+(d{y_{i}})^{2})
\end{array}
\]
and some compatible tensor $\tilde{\Phi}$. Such space $(\mathbb{R}^{5}%
,\tilde{\Phi},\xi,\tilde{\eta},\tilde{g})$ is usually referred to as the standard
Sasakian space form with $\mathcal{H}(X)=-3$.
\end{remark}

The choice of our $\eta$-form is actually made for the convenience of
calculation due to its symmetry. Moreover, the next lemma combined with
proposition \ref{curvature} states that our ambient space $(\mathbb{R}%
^{5},\Phi,\xi,\eta,g)$ is congruent to the standard Sasakian space form
$(\mathbb{R}^{5},\tilde{\Phi},\xi,\tilde{\eta},\tilde{g})$.

\begin{lemma}
\label{prop2.1} $(\mathbb{R}^{5},\Phi,\xi,\eta,g)$ is a Sasakian manifold with
$\mathcal{H}(X)=-3$ for any vector $X$ orthogonal to the Reeb vector field
$\xi$.
\end{lemma}

\begin{proof}
Firstly, we construct an adapted frame field $\{E_{i}\}_{i=1}^{5}$ explicitly
as%
\begin{equation}%
\begin{array}
[c]{c}%
E_{1}=2\frac{\partial}{\partial x_{1}}+y_{1}\frac{\partial}{\partial z},\text{
}E_{2}=2\frac{\partial}{\partial x_{2}}+y_{2}\frac{\partial}{\partial z},\\
E_{3}=-2\frac{\partial}{\partial y_{1}}+x_{1}\frac{\partial}{\partial
z},\text{ }E_{4}=-2\frac{\partial}{\partial y_{2}}+x_{2}\frac{\partial
}{\partial z},\text{ }E_{5}=\xi=2\frac{\partial}{\partial z},
\end{array}
\label{eq2.2}%
\end{equation}
which satisfy
\begin{equation}%
\begin{array}
[c]{c}%
\eta\left(  E_{i}\right)  =\eta\left(  E_{2+i}\right)  =0,\\
\Phi\left(  E_{i}\right)  =E_{2+i},\text{ \textrm{for} }i=1,\text{ }2,\\
\Phi(\xi)=0,\\
g(E_{i},E_{j})=\delta_{ij},
\end{array}
\label{eq2.3}%
\end{equation}
such frame is called a $\Phi$-basis of $(\mathbb{R}^{5},\Phi,\xi,\eta,g)$.
Utilizing this frame, one can verify that it satisfies the Sasakian condition
(2.2)-(2.5) one by one.

For (\ref{eq2.2a}), it is easy to see
\[
\Phi^{2}(\cdot,\xi)+I(\cdot,\xi)=\xi\otimes\eta(\cdot,\xi)
\]
holds; on the other hand
\[
\eta\left(  E_{i}\right)  =\eta\left(  E_{n+i}\right)  =0,\text{ \textrm{for}
}i=1,\text{ }2,
\]
it follows
\[%
\begin{array}
[c]{l}%
\Phi^{2}\left(  E_{i}\right)  =\Phi\left(  E_{2+i}\right)  =\Phi\left(
-2\frac{\partial}{\partial y_{i}}+x_{i}\frac{\partial}{\partial
z}\right) \\
=-2\Phi(\frac{\partial}{\partial y_{i}})=-2\frac{\partial}{\partial x_{i}%
}-y_{i}\frac{\partial}{\partial z}=-E_{i},\text{ }\mathrm{for}\text{ }i=1,2.
\end{array}
\]

For (\ref{eq2.3a}),%
\[%
\begin{array}
[c]{c}%
d\eta=\frac{1}{2}\sum_{i=1}^{2}dx_{i}\wedge{dy_{i}}=\frac{1}{2}(dx_{1}%
\wedge{dy_{1}}+dx_{2}\wedge{dy_{2}}),
\end{array}
\]
so%
\[%
\begin{array}
[c]{c}%
\frac{1}{2}d\eta(\xi,E_{i})=0=g(\xi,\Phi(E_{i})),
\end{array}
\]%
\[%
\begin{array}
[c]{c}%
\frac{1}{2}d\eta(\cdot,\xi)=0=g(\cdot,\Phi(\xi)),
\end{array}
\]
and%
\[%
\begin{array}
[c]{c}%
\frac{1}{2}d\eta(\frac{\partial}{\partial x_{1}},\frac{\partial}{\partial
x_{2}})=0=g(\frac{\partial}{\partial x_{1}},\Phi(\frac{\partial}{\partial
x_{2}}))=-g(\frac{\partial}{\partial x_{1}},\frac{\partial}{\partial y_{1}%
})+x_{2}g(\frac{\partial}{\partial x_{1}},\frac{\partial}{\partial z}).
\end{array}
\]
Similarly, we have
\[%
\begin{array}
[c]{c}%
\frac{1}{2}d\eta(\frac{\partial}{\partial y_{1}},\frac{\partial}{\partial
y_{2}})=g(\frac{\partial}{\partial y_{1}},\Phi(\frac{\partial}{\partial y_{2}%
}))=0,
\end{array}
\]%
\[%
\begin{array}
[c]{c}%
\frac{1}{2}d\eta(\frac{\partial}{\partial x_{1}},\frac{\partial}{\partial
y_{1}})=\frac{1}{4}=g(\frac{\partial}{\partial x_{1}},\Phi(\frac{\partial
}{\partial y_{1}}))=g(\frac{\partial}{\partial x_{1}},\frac{\partial}{\partial
x_{1}})+\frac{y_{1}}{2}g(\frac{\partial}{\partial x_{1}},\frac{\partial
}{\partial z}),
\end{array}
\]
Also for $\frac{1}{2}d\eta(\frac{\partial}{\partial x_{i}},\frac{\partial
}{\partial y_{i}})$, we can obtain $\frac{1}{2}d\eta(\cdot,\cdot)=g(\cdot
,\Phi(\cdot))$.

For (\ref{eq2.4a}),%
\[
g(\xi,\xi)=1=g(\Phi(\xi),\Phi(\xi))+\eta(\xi)\otimes\eta(\xi),
\]%
\[
g(E_{i},\xi)=0=g(\Phi(E_{i}),\Phi(\xi))+\eta(E_{i})\otimes\eta(\xi),
\]
and
\begin{align*}
g(E_{1},E_{1})  &  =1=g(\Phi(E_{1}),\Phi(E_{1}))+\eta(E_{1})\otimes\eta
(E_{1}),\\
g(E_{3},E_{3})  &  =1=g(\Phi(E_{3}),\Phi(E_{3}))+\eta(E_{3})\otimes\eta
(E_{3}),\\
g(E_{1},E_{2})  &  =0=g(\Phi(E_{1}),\Phi(E_{2}))+\eta(E_{1})\otimes\eta
(E_{2}),\\
g(E_{3},E_{4})  &  =0=g(\Phi(E_{3}),\Phi(E_{4}))+\eta(E_{3})\otimes\eta
(E_{4}),
\end{align*}
Similar for other $g(E_{i},E_{i})$ yields $g(\cdot,\cdot)=g(\Phi\cdot
,\Phi\cdot)+\eta\otimes\eta.$

For (\ref{eq2.5a}),%
\[%
\begin{array}
[c]{lll}%
d\eta\otimes\xi & = & (dx_{1}\wedge{dy_{1}}+dx_{2}\wedge{dy_{2}})\otimes
\frac{\partial}{\partial z}\\
& = & (dx_{1}\otimes{dy_{1}}+dx_{2}\otimes{dy_{2}}-dy_{1}\otimes{dx_{1}%
}-dy_{2}\otimes{dx_{2}})\otimes\frac{\partial}{\partial z}.
\end{array}
\]
By the definition of tensor $\Phi$ in (\ref{eq2.1}), we have%
\[%
\begin{array}
[c]{c}%
\Phi(\frac{\partial}{\partial x_{1}})=-\frac{\partial}{\partial y_{1}}%
+\frac{x_{1}}{2}\frac{\partial}{\partial z},\text{ \ }\Phi(\frac{\partial
}{\partial x_{2}})=-\frac{\partial}{\partial y_{2}}+\frac{x_{2}}{2}%
\frac{\partial}{\partial z},\\
\Phi(\frac{\partial}{\partial y_{1}})=\frac{\partial}{\partial x_{1}}%
+\frac{y_{1}}{2}\frac{\partial}{\partial z},\text{ \ }\Phi(\frac{\partial
}{\partial y_{2}})=\frac{\partial}{\partial x_{2}}+\frac{y_{2}}{2}%
\frac{\partial}{\partial z},\text{ \ }\Phi(\frac{\partial}{\partial z})=0,
\end{array}
\]
and%
\[%
\begin{array}
[c]{lll}%
\Phi & = & dx^{1}\otimes(-\frac{\partial}{\partial y_{1}}+\frac{x_{1}}{2}%
\frac{\partial}{\partial z})+dx^{2}\otimes(-\frac{\partial}{\partial y_{2}%
}+\frac{x_{2}}{2}\frac{\partial}{\partial z})\\
&  & +dy^{1}\otimes(\frac{\partial}{\partial x_{1}}+\frac{y_{1}}{2}%
\frac{\partial}{\partial z})+dy^{2}\otimes(\frac{\partial}{\partial x_{2}%
}+\frac{y_{2}}{2}\frac{\partial}{\partial z})\\
& = & -dx^{1}\otimes\frac{\partial}{\partial y_{1}}-dx^{2}\otimes
\frac{\partial}{\partial y_{2}}+dy^{1}\otimes\frac{\partial}{\partial x_{1}%
}+dy^{2}\otimes\frac{\partial}{\partial x_{2}}\\
&  & +(\frac{x_{1}}{2}dx^{1}+\frac{x_{2}}{2}dx^{2}+\frac{y_{1}}{2}dy^{1}%
+\frac{y_{2}}{2}dy^{2})\otimes\frac{\partial}{\partial z}\\
& = & J_{\mathbb{C}^{2}}+T_{0},
\end{array}
\]
where $T_{0}$ is the part parallel to $\xi$.

It follows that
\[
\mathcal{N}_{\Phi}(X,Y)=\Phi^{2}[X,Y]+[\Phi{X},\Phi{Y}]-\Phi\lbrack\Phi
{X},Y]-\Phi\lbrack X,\Phi{Y}],
\]
$\mathcal{N}_{\Phi}=0$ when $X=Y$ and $\mathcal{N}_{\Phi}(X,Y)=-\mathcal{N}%
_{\Phi}(Y,X)$, it implies that%
\[%
\begin{array}
[c]{lll}%
\mathcal{N}_{\Phi}(\frac{\partial}{\partial x_{1}},\frac{\partial}{\partial
x_{2}}) & = & \Phi^{2}[\frac{\partial}{\partial x_{1}},\frac{\partial
}{\partial x_{2}}]+[\Phi\frac{\partial}{\partial x_{1}},\Phi\frac{\partial
}{\partial x_{2}}]-\Phi\lbrack\Phi\frac{\partial}{\partial x_{1}}%
,\frac{\partial}{\partial x_{2}}]-\Phi\lbrack\frac{\partial}{\partial x_{1}%
},\Phi\frac{\partial}{\partial x_{2}}]\\
& = & 0+[-\frac{\partial}{\partial y_{1}}+\frac{x_{1}}{2}\frac{\partial
}{\partial z},-\frac{\partial}{\partial y_{2}}+\frac{x_{2}}{2}\frac{\partial
}{\partial z}]-\Phi\lbrack-\frac{\partial}{\partial y_{1}}+\frac{x_{1}}%
{2}\frac{\partial}{\partial z},\frac{\partial}{\partial x_{2}}]\\
&  & -\Phi\lbrack\frac{\partial}{\partial x_{1}},-\frac{\partial}{\partial
y_{2}}+\frac{x_{2}}{2}\frac{\partial}{\partial z}]\\
& = & [-\frac{\partial}{\partial y_{1}},\frac{x_{2}}{2}\frac{\partial
}{\partial z}]+[\frac{x_{1}}{2}\frac{\partial}{\partial z},-\frac{\partial
}{\partial y_{2}}]+[\frac{x_{1}}{2}\frac{\partial}{\partial z},\frac{x_{2}}%
{2}\frac{\partial}{\partial z}]\\
&  & -\Phi\lbrack\frac{x_{1}}{2}\frac{\partial}{\partial z},\frac{\partial
}{\partial x_{2}}]-\Phi\lbrack\frac{\partial}{\partial x_{1}},\frac{x_{2}}%
{2}\frac{\partial}{\partial z}]\\
& = & 0,
\end{array}
\]
and
\[%
\begin{array}
[c]{lll}%
\mathcal{N}_{\Phi}(\frac{\partial}{\partial x_{1}},\frac{\partial}{\partial
y_{1}}) & = & \Phi^{2}[\frac{\partial}{\partial x_{1}},\frac{\partial
}{\partial y_{1}}]+[\Phi\frac{\partial}{\partial x_{1}},\Phi\frac{\partial
}{\partial y_{1}}]-\Phi\lbrack\Phi\frac{\partial}{\partial x_{1}}%
,\frac{\partial}{\partial y_{1}}]-\Phi\lbrack\frac{\partial}{\partial x_{1}%
},\Phi\frac{\partial}{\partial y_{1}}]\\
& = & 0+[-\frac{\partial}{\partial y_{1}}+\frac{x_{1}}{2}\frac{\partial
}{\partial z},\frac{\partial}{\partial x_{1}}+\frac{y_{1}}{2}\frac{\partial
}{\partial z}]-\Phi\lbrack-\frac{\partial}{\partial y_{1}}+\frac{x_{1}}%
{2}\frac{\partial}{\partial z},\frac{\partial}{\partial y_{1}}]\\
&  & -\Phi\lbrack\frac{\partial}{\partial x_{1}},-\frac{\partial}{\partial
x_{1}}+\frac{y_{1}}{2}\frac{\partial}{\partial z}]\\
& = & [-\frac{\partial}{\partial y_{1}},\frac{y_{1}}{2}\frac{\partial
}{\partial z}]+[\frac{x_{1}}{2}\frac{\partial}{\partial z},\frac{\partial
}{\partial x_{1}}]+[\frac{x_{1}}{2}\frac{\partial}{\partial z},\frac{y_{1}}%
{2}\frac{\partial}{\partial z}]\\
&  & -\Phi\lbrack\frac{x_{1}}{2}\frac{\partial}{\partial z},\frac{\partial
}{\partial y_{1}}]-\Phi\lbrack\frac{\partial}{\partial x_{1}},\frac{y_{1}}%
{2}\frac{\partial}{\partial z}]\\
& = & 0,
\end{array}
\]
thus
\[%
\begin{array}
[c]{c}%
\mathcal{N}_{\Phi}=-(dx_{1}\wedge{dy_{1}}+dx_{2}\wedge{dy_{2}})\otimes
\frac{\partial}{\partial z}=-d\eta\otimes\xi.
\end{array}
\]
Hence, $(\mathbb{R}^{5},\eta,\xi,\Phi,g)$ is a Sasakian manifold.

Next, we apply Koszul's formula to calculate Riemannian connections of frames:%
\begin{equation}%
\begin{array}
[c]{lll}%
2g(\nabla_{X}Y,Z) & = & Xg(Y,Z)+Yg(X,Z)-Zg(X,Y)\\
&  & +g([X,Y],Z)+g([Z,X],Y)-g([Y,Z],X),
\end{array}
\label{eq2.4}%
\end{equation}
firstly, we calculate Lie brackets:
\[
\lbrack E_{1},E_{2}]=[E_{1},E_{4}]=[E_{2},E_{3}]=[E_{3},E_{4}]=0,
\]
and%
\[%
\begin{array}
[c]{lll}%
\lbrack E_{1},E_{3}] & = & [2{\partial x_{1}}+y_{1}{\partial z},-2{\partial
y_{1}}+x_{1}{\partial z}]\\
& = & -4[{\partial x_{1}},{\partial y_{1}}]+[y_{1}{\partial z},x_{1}{\partial
z}]+2[{\partial x_{1}},x_{1}{\partial z}]-2[y_{1}{\partial z},{\partial y_{1}%
}]\\
& = & 4{\partial z}=2\xi,
\end{array}
\]
also $[E_{2},E_{4}]=2\xi,$ $[\xi,\xi]=0$, $[E_{i},E_{i}]=0,$ and $[E_{1}%
,\xi]=0$ for $i=1,\cdots,4$.

Then we obtain by (\ref{eq2.4})
\begin{align*}
\overline{\nabla}_{E_{1}}E_{1}  &  =\overline{\nabla}_{E_{2}}E_{2}%
=\overline{\nabla}_{E_{3}}E_{3}=\overline{\nabla}_{E_{4}}E_{4}=\overline
{\nabla}_{\xi}\xi=0,\\
\overline{\nabla}_{E_{1}}E_{2}  &  =\overline{\nabla}{E_{2}}E_{1}=0,\text{
}\overline{\nabla}_{E_{1}}E_{3}=\xi=-\overline{\nabla}_{E_{3}}E_{1},\\
\overline{\nabla}_{E_{1}}E_{4}  &  =\overline{\nabla}_{E_{4}}E_{1}=0,\text{
}\overline{\nabla}_{E_{1}}\xi=-E_{3}=-\overline{\nabla}_{\xi}E_{1},\\
\overline{\nabla}_{E_{2}}E_{3}  &  =\overline{\nabla}_{E_{3}}E_{2}=0,\text{
}\overline{\nabla}_{E_{2}}E_{4}=\xi=-\overline{\nabla}_{E_{4}}E_{2},\\
\overline{\nabla}_{E_{2}}\xi &  =-E_{4}=-\overline{\nabla}_{\xi}E_{2},\text{
}\overline{\nabla}_{E_{3}}E_{4}=\overline{\nabla}_{E_{4}}E_{3}=0,\\
\overline{\nabla}_{E_{3}}\xi &  =E_{1}=-\overline{\nabla}_{\xi}E_{3},\text{
}\overline{\nabla}_{E_{4}}\xi=E_{2}=-\overline{\nabla}_{\xi}E_{4}.
\end{align*}

Now straightforward calculation will show that $\Phi$-sectional curvature of
$(\mathbb{R}^{5},\Phi,\xi,\eta,g)$ is $-3$.
\end{proof}

\begin{definition}
A submanifold $F:L^{n}\rightarrow M^{2n+1}$ of a Sasakian manifold
$(M^{2n+1},\Phi,\xi,\eta,g)$ is called isotropic if
\[
\eta|_{TL}=0\text{ \textrm{or} }F^{\ast}\eta=0.
\]
In particular, $F^{\ast}d\eta=0$ as well. A Legendrian submanifold $L$ is a
maximally isotropic submanifold of dimension $n.$
\end{definition}

The Legendrian manifold is analogous to the Lagrangian manifold in many ways.
Now we introduce the so-called Lagrangian projection:
\[
\Pi:(\mathbb{R}^{2n+1},\Phi,\xi,\eta,g)\rightarrow(\mathbb{C}^{n}%
,J,\omega_{\mathrm{can}})
\]
as follows:
\[%
\begin{array}
[c]{ccc}
&  & ((\mathbb{R}^{2n+1},\Phi,\xi,\eta,g))\ni(x_{1},y_{1},\cdots,x_{n}%
,y_{n},z)\\
& F\nearrow\curvearrowright & \downarrow\Pi\text{
\ \ \ \ \ \ \ \ \ \ \ \ \ \ \ \ \ \ \ \ \ \ \ }\\
L^{n} & \overset{f}{\longrightarrow} & (\mathbb{C}^{n},J,\omega_{\mathrm{can}%
})\ni(x_{1},y_{1},\cdots,x_{n},y_{n}).
\end{array}
\]
Then, since $\Pi^{\ast}(\omega_{\mathrm{can}})=\frac{1}{2}d\eta$ and
$\Pi_{\ast}(\xi)=0$%
\[
f:L^{n}\rightarrow(\mathbb{C}^{n},J,\omega_{\mathrm{can}})
\]
is a Lagrangian submanifold with $f=\Pi\circ F$.

Note that if $F_{t}:L^{n}\times\lbrack0,T)\rightarrow\mathbb{R}^{2n+1}$ is the
maximal existence Legendrian mean curvature flow
\[%
\begin{array}
[c]{c}%
\frac{\partial F_{t}}{\partial t}=H_{\mathbb{R}^{2n+1}}-2\alpha\xi,
\end{array}
\]
then $f_{t}=\Pi\circ F_{t}:L^{n}\times\lbrack0,T)\rightarrow\mathbb{C}^{n}$ is
the solution of the Lagrangian mean curvature flow
\[%
\begin{array}
[c]{c}%
\frac{\partial f_{t}}{\partial t}=\frac{\partial(\Pi\circ F_{t})}{\partial
t}=\Pi_{\ast}(\frac{\partial F_{t}}{\partial t})=H_{\mathbb{C}^{n}}.
\end{array}
\]
Here the Legendre angle $\alpha$ of $F_{t}$ is the one defined in
(\ref{angle}) by $F_{t}^{\ast}\Omega^{T}=e^{-i\alpha}Vol_{L_{t}}$, where
\[
\Omega^{T}=\frac{1}{4}dz_{1}\wedge dz_{2}%
\]
is the transverse holomorphic volume form, normalized so that $|\Omega
^{T}|_{g^{T}}=1$; equivalently $\omega^{T}\lrcorner H=d\alpha$, i.e.
$H=\Phi\nabla\alpha$, where $\omega^{T}=g^{T}(\cdot,\Phi\cdot)=\frac{1}%
{2}d\eta$ is the transverse K\"{a}hler form. Observe also that $\Pi_{\ast
}(\alpha\xi)=0$: the vertical term of (\ref{c}) is invisible after the
Lagrangian projection, and the coefficient in front of it can only be detected
through the Legendrian condition itself.
Conversely, one can recover a Legendrian manifold from a Lagrangian manifold
via the lifting theorem:

\begin{proposition}
[Theorem 2.1 of \cite{bmv}]\label{P51} Given a Lagrangian immersion
$f:\Sigma^{n}\rightarrow(\mathbb{C}^{n},J,\omega_{\mathrm{can}}):$%
\[%
\begin{array}
[c]{ccc}
&  & (\mathbb{C}^{n}\times\mathbb{S}^{1},\Phi,\xi,\eta,g)\ni(x_{1}%
,y_{1},\cdots,x_{n},y_{n},e^{iz})\\
& \bar{F}\nearrow\curvearrowright & \downarrow\bar{\Pi}\text{
\ \ \ \ \ \ \ \ \ \ \ \ \ \ \ \ \ \ \ \ \ \ }\\
L^{n} & \overset{f}{\longrightarrow} & (\mathbb{C}^{n},J,\omega_{\mathrm{can}%
})\ni(x_{1},y_{1},\cdots,x_{n},y_{n}).
\end{array}
\]
Then

\begin{enumerate}
\item $\bar{F}$ is a local Legendrian immersion in $\mathbb{C}^{n}%
\times\mathbb{S}^{1}$ with
\begin{equation}%
\begin{array}
[c]{c}%
z(x):=\int_{\gamma}\eta,
\end{array}
\label{a}%
\end{equation}
where $(\gamma(x),e^{iz(x)})\in L^{n}\times\mathbb{S}^{1}$ and $z(x_{0})=0.$

\item $\bar{F}$ is a global Legendrian immersion if
\begin{equation}%
\begin{array}
[c]{c}%
z(x)=\int_{\gamma}\eta\in2\pi\mathbb{Z},
\end{array}
\label{b}%
\end{equation}
for all $[\gamma]\in H_{1}(L,\mathbb{Z}).$

\item $\bar{F}$ is an embedding if
\[%
\begin{array}
[c]{c}%
\int_{\gamma}\eta\neq0,\text{\ }\mathrm{(mod}\text{ }2\pi\mathrm{)}%
\end{array}
\]
along any path joining a pair of double points in the sense that for all
points $x_{0},$ $x_{1}\in\Sigma$ such that $f(x_{0})=f(x_{1})$ and any path
$\gamma$ from $x_{0}$ to $x_{1}$ in $L$.
\end{enumerate}
\end{proposition}

\begin{remark}
We call $\bar{F}$ the corresponding Legendrian immersion of $F$.
\end{remark}

\section{Examples of Legendrian self-similar solution}

In this section, we will construct examples for Legendrian self-similar
solution. Our method is applying the lifting theorem \ref{P51} and known
models constructed by Anciaux \cite{a} and Joyce-Lee-Tsui \cite{jlt} for the
Lagrangian case. Here we split the proof of Theorem \ref{Ex} into 2 parts.

Let $\lambda_{1},\lambda_{2},C\in\mathbb{R} \backslash\{0\}$, $I$ be an open
interval in $\mathbb{R} $, and $\alpha:I\rightarrow\mathbb{R} $ or
$\alpha:I\rightarrow\mathbb{R}{/ 2 \pi\mathbb{Z}} $ and $w_{1}, w_{2}:
I \rightarrow\mathbb{C} \backslash\{0\} $ be smooth functions.

\begin{lemma}
Write $w_{j}=r_{j}(s)e^{i\phi_{j}(s)},r_{j}=|w_{j}|$ and $\phi=\phi_{1}%
+\phi_{2}$, for functions $r_{j}:I\rightarrow(0,\infty)$ and $\phi_{1}%
,\phi_{2},\phi:I\rightarrow\mathbb{R}$ or $\mathbb{R}{/2}\pi\mathbb{Z}$. Fix
$s_{0}\in{I}$. Define $u:I\rightarrow\mathbb{R}$ by
\begin{equation}%
\begin{array}
[c]{c}%
u(s)=2\int_{s_{0}}^{s}r_{1}(t)r_{2}(t)\cos(\alpha(t)+\phi(t))dt.
\end{array}
\label{eq15}%
\end{equation}
For the above $u$, we have $r_{j}^{2}(s)=\rho_{j}+\lambda_{j}u(s),$ $j=1,2,$
$s\in I,$ $\rho_{j}=r_{j}^{2}\left(  s_{0}\right)  $. Define a degree 2
polynomial $Q(u)$ by $Q(u)=\prod_{j=1}^{2}\left(  \rho_{j}+\lambda
_{j}u\right)  $. Suppose that%
\begin{equation}
\left\{
\begin{array}
[c]{l}%
\frac{dw_{1}}{ds}=\lambda_{1}e^{-i\alpha(s)}\overline{w}_{2},\text{
\ }\frac{dw_{2}}{ds}=\lambda_{2}e^{-i\alpha(s)}\overline{w}_{1},\\
\frac{d\phi_{j}}{ds}=-\frac{\lambda_{j}}{\rho_{j}+\lambda_{j}u}%
Q(u)^{1/2}\sin(\alpha+\phi),\text{ }j=1,2,\\
\frac{dz}{ds}=\frac{C}{2}Q(u)^{1/2}\sin(\alpha+\phi),
\end{array}
\right.  \label{eq16}%
\end{equation}
hold in $I$. Then the submanifold $\tilde{L}$ in $\mathbb{C}{^{2}}{{\times
}\mathbb{R}}$ given by
\begin{equation}%
\begin{array}
[c]{c}%
\tilde{L}=\left\{  \left(  x_{1}w_{1}(s),x_{2}w_{2}(s),z(s)\right)  :s\in
I,\text{ }x_{j}\in\mathbb{R},\text{ }\sum_{j=1}^{2}\lambda_{j}x_{j}%
^{2}=C\right\}
\end{array}
\label{eq17}%
\end{equation}
is Legendrian.
\end{lemma}

\begin{proof}
Using equations (\ref{eq16}), $w_{j}(s)=r_{j}(s)e^{i\phi_{j}(s)},$
$r_{j}(s)=|w_{j}(s)|$ and $|w_{j}|^{2}=w_{j}\cdot\bar{w}_{j}$ for $j=1,2$ we
have%
\[%
\begin{array}
[c]{lll}%
\frac{d}{ds}r_{1}^{2} & = & \frac{d}{ds}|w_{1}|^{2}=\frac{d}{ds}(w_{1}%
\cdot\bar{w}_{1})=w_{1}\frac{d}{ds}\bar{w}_{1}+\overline{w}_{1}\frac{d}%
{ds}w_{1}\\
& = & \lambda_{1}(w_{1}e^{i\alpha}w_{2}+\overline{w}%
_{1}e^{-i\alpha}\overline{w}_{2})=\lambda_{1}r_{1}r_{2}%
(e^{i\alpha}e^{i\phi}+e^{-i\alpha}e^{-i\phi})\\
& = & \lambda_{1}r_{1}r_{2}(e^{i(\alpha+\phi)}+e^{-i(\alpha+\phi)}%
)=2\lambda_{1}r_{1}r_{2}\cos(\alpha+\phi)\\
& = & \lambda_{1}\frac{du}{ds},
\end{array}
\]
Similarly,
\[%
\begin{array}
[c]{c}%
\frac{d}{ds}r_{2}^{2}=\lambda_{2}\frac{du}{ds}.
\end{array}
\]
Then $\frac{d}{ds}(r_{j}^{2}-\lambda_{j}u)=0$. Hence $r_{j}^{2}-\lambda_{j}u$
is constant in $I$, and at $s=s_{0}$, we have $u(s_{0})=0$ and $r_{j}%
^{2}(s_{0})=\rho_{j}$. So $r_{j}^{2}(s)=\rho_{j}+\lambda_{j}u(s)$ for
$j=1,2,$ $s\in{I},$ and the definition of $Q(u)$ imply that $Q(u)=r_{1}%
^{2}r_{2}^{2}.$

To facilitate matters, we can express $\tilde{L}$ in the following form with
$w_{j}=r_{j}(s)e^{i\phi_{j}(s)}$:
\begin{equation}
\tilde{L}=\left(  x_{1}r_{1}(s)\cos\phi_{1}(s),x_{1}r_{1}(s)\sin\phi
_{1}(s),x_{2}r_{2}(s)\cos\phi_{2}(s),x_{2}r_{2}(s)\sin\phi_{2}%
(s),z(t,s)\right)  , \label{eq18}%
\end{equation}
where the ambient space $\mathbb{R}^{5}$ is using coordinates $(\bar{x_{1}%
},\bar{y_{1}},\bar{x_{2}},\bar{y_{2}},\bar{z})$. From (\ref{eq18}) and the
third equation of (\ref{eq16}), it follows that%
\[%
\begin{array}
[c]{lll}%
2dz & = & x_{1}r_{1}\sin\phi_{1}d(x_{1}r_{1}\cos\phi_{1})-x_{1}r_{1}\cos
\phi_{1}d(x_{1}r_{1}\sin\phi_{1})\\
&  & +x_{2}r_{2}\sin\phi_{2}d(x_{2}r_{2}\cos\phi_{2})-x_{2}r_{2}\cos\phi
_{2}d(x_{2}r_{2}\sin\phi_{2})\\
& = & x_{1}r_{1}\sin\phi_{1}(r_{1}\cos\phi_{1}dx_{1}+x_{1}\cos\phi_{1}%
dr_{1}-x_{1}r_{1}\sin\phi_{1}d\phi_{1})\\
&  & -x_{1}r_{1}\cos\phi_{1}(r_{1}\sin\phi_{1}dx_{1}+x_{1}\sin\phi_{1}%
dr_{1}+x_{1}r_{1}\cos\phi_{1}d\phi_{1})\\
&  & +x_{2}r_{2}\sin\phi_{2}(r_{2}\cos\phi_{2}dx_{2}+x_{2}\cos\phi_{2}%
dr_{2}-x_{2}r_{2}\sin\phi_{2}d\phi_{2})\\
&  & -x_{2}r_{2}\cos\phi_{2}(r_{2}\sin\phi_{2}dx_{2}+x_{2}\sin\phi_{2}%
dr_{2}+x_{2}r_{2}\cos\phi_{2}d\phi_{2})\\
& = & -x_{1}^{2}r_{1}^{2}d\phi_{1}-x_{2}^{2}r_{2}^{2}d\phi_{2}=(\lambda
_{1}x_{1}^{2}+\lambda_{2}x_{2}^{2})Q(u)^{1/2}\sin(\alpha+\phi)ds\\
& = & CQ(u)^{1/2}\sin(\alpha+\phi).
\end{array}
\]
On the other hand, $2dz=2z_{s}ds+2z_{t}dt=CQ(u)^{1/2}\sin(\alpha+\phi)$ follows from the last equation of (\ref{eq16}), thus $\tilde{L}^{\ast}%
\eta=0$, and it can be deduced that $z_{t}=0$. Hence $\tilde{L}$ is a
Legendrian submanifold in $\mathbb{C}{^{2}}{{\times}\mathbb{R}}$.
\end{proof}

\begin{theorem}
Under the same condition of the above lemma, if it holds
\begin{equation}%
\begin{array}
[c]{c}%
\frac{d\alpha}{ds}=-\frac{a}{4}Cr_{1}r_{2}\sin(\alpha+\phi)
\end{array}
\label{eq19}%
\end{equation}
in $I$, then $\tilde{L}$ satisfies (\ref{CCC}) and is a self-similar solution
to the Legendrian mean curvature flow. Moreover, for some $A\in\mathbb{R}$ the
equations (\ref{eq16}) and (\ref{eq19}) admit the first integral
\begin{equation}
Q(u)^{1/2}e^{a{Cu}/8}\sin(\alpha+\phi)=A. \label{eq20}%
\end{equation}

\end{theorem}

\begin{proof}
We only need to discuss the following two cases of $\lambda_{1},$ $\lambda
_{2}$ and $C$ to cover all possibilities:

i) $\lambda_{1}>0,$ $\lambda_{2}>0,$ $C>0$; ii) $\lambda_{1}>0,$ $\lambda
_{2}<0,$ $C>0$.

Now let us first discuss i). We choose the orthonormal frame $\{E_{i}%
,\xi\}_{i=1}^{4}$. Using $w_{j}=r_{j}(s)e^{i\phi_{j}(s)}$ and parameterize the
equation $\sum_{j=1}^{2}\lambda_{j}x_{j}^{2}=C$, for $j=1,2$. Under natural
coordinates, we have%
\begin{equation}%
\begin{array}
[c]{lll}%
\tilde{L} & = & \left(  \sqrt{\frac{C}{\lambda_{1}}}r_{1}(s)\cos t\cos\phi
_{1}(s),\sqrt{\frac{C}{\lambda_{1}}}r_{1}(s)\cos t\sin\phi_{1}(s),\right. \\
&  & \left.  \sqrt{\frac{C}{\lambda_{2}}}r_{2}(s)\sin t\cos\phi_{2}%
(s),\sqrt{\frac{C}{\lambda_{2}}}r_{2}(s)\sin t\sin\phi_{2}(s),z(t,s)\right)  ,
\end{array}
\label{eq21}%
\end{equation}
with adapted frame (\ref{eq2.2}), $\tilde{L}$ can be expressed as%
\begin{equation}%
\begin{array}
[c]{lll}%
\tilde{L} & = & (\sqrt{\frac{C}{\lambda_{1}}}r_{1}\cos t\cos\phi_{1}%
,\sqrt{\frac{C}{\lambda_{1}}}r_{1}\cos t\sin\phi_{1},\\
&  & \sqrt{\frac{C}{\lambda_{2}}}r_{2}\sin t\cos\phi_{2},\sqrt{\frac
{C}{\lambda_{2}}}r_{2}\sin t\sin\phi_{2},z)\\
& = & \sqrt{\frac{C}{\lambda_{1}}}r_{1}\cos t\cos\phi_{1}(\frac{1}{2}%
E_{1}-\frac{\sqrt{C}}{4\sqrt{\lambda_{1}}}r_{1}\cos t\sin\phi_{1}\xi)\\
&  & +\sqrt{\frac{C}{\lambda_{1}}}r_{1}\cos t\sin\phi_{1}(-\frac{1}{2}%
E_{3}+\frac{\sqrt{C}}{4\sqrt{\lambda_{1}}}r_{1}\cos t\sin\phi_{1}\xi)\\
&  & +\sqrt{\frac{C}{\lambda_{2}}}r_{2}\sin t\cos\phi_{2}(\frac{1}{2}%
E_{2}-\frac{\sqrt{C}}{4\sqrt{\lambda_{1}}}r_{2}\sin t\sin\phi_{2}\xi)\\
&  & +\sqrt{\frac{C}{\lambda_{2}}}r_{2}\sin t\sin\phi_{2}(-\frac{1}{2}%
E_{4}+\frac{\sqrt{C}}{4\sqrt{\lambda_{1}}}r_{2}\sin t\cos\phi_{2}\xi)+\frac
{1}{2}z\xi\\
& = & \frac{r_{1}}{2}\sqrt{\frac{C}{\lambda_{1}}}\cos t\cos\phi_{1}E_{1}%
+\frac{r_{2}}{2}\sqrt{\frac{C}{\lambda_{2}}}\sin t\cos\phi_{2}E_{2}\\
&  & -\frac{r_{1}}{2}\sqrt{\frac{C}{\lambda_{1}}}\cos t\sin\phi_{1}E_{3}%
-\frac{r_{2}}{2}\sqrt{\frac{C}{\lambda_{2}}}\sin t\sin\phi_{2}E_{4}+\frac
{1}{2}z\xi.
\end{array}
\label{eq22}%
\end{equation}
Using the condition (\ref{eq16}), $w_{j}(s)=r_{j}(s)e^{i\phi_{j}(s)}%
,r_{j}(s)=|w_{j}(s)|$ and $\phi=\phi_{1}+\phi_{2}$, we obtain
\[%
\begin{array}
[c]{lll}%
\frac{d}{ds}r_{1} & = & \frac{d}{ds}|w_{1}|=\frac{d}{ds}\sqrt{\langle
w_{1},\bar{w_{1}}\rangle}=\frac{1}{2r_{1}}(\langle\frac{dw_{1}}{ds},\bar
{w_{1}}\rangle+\langle w_{1},\frac{d\bar{w_{1}}}{ds}\rangle)\\
& = & \frac{\lambda_{1}r_{2}}{2}(e^{i(\alpha+\phi)}+e^{-i(\alpha+\phi)}%
)=\lambda_{1}r_{2}\cos(\alpha+\phi),
\end{array}
\]
and also $\frac{dr_{2}}{ds}=\lambda_{2}r_{1}\cos(\alpha+\phi)$. Thus
the tangent vectors are%
\begin{equation}%
\begin{array}
[c]{lll}%
\frac{\partial\tilde{L}}{\partial t} & = & (-\sqrt{\frac{C}{\lambda_{1}}}%
r_{1}\sin t\cos\phi_{1},-\sqrt{\frac{C}{\lambda_{1}}}r_{1}\sin t\sin\phi
_{1},\\
&  & \sqrt{\frac{C}{\lambda_{2}}}r_{2}\cos t\cos\phi_{2},\sqrt{\frac
{C}{\lambda_{2}}}r_{2}\cos t\sin\phi_{2},0)\\
& = & -\frac{r_{1}}{2}\sqrt{\frac{C}{\lambda_{1}}}\sin t\cos\phi_{1}%
E_{1}+\frac{r_{2}}{2}\sqrt{\frac{C}{\lambda_{2}}}\cos t\cos\phi_{2}E_{2}\\
&  & +\frac{r_{1}}{2}\sqrt{\frac{C}{\lambda_{1}}}\sin t\sin\phi_{1}E_{3}%
-\frac{r_{2}}{2}\sqrt{\frac{C}{\lambda_{2}}}\cos t\sin\phi_{2}E_{4},
\end{array}
\label{eq23}%
\end{equation}
and%
\begin{equation}%
\begin{array}
[c]{lll}%
\frac{\partial\tilde{L}}{\partial s} & = & (\sqrt{\frac{C}{\lambda_{1}}%
}\lambda_{1}r_{2}\cos t[\cos\phi_{1}\cos(\alpha+\phi)+\sin\phi_{1}\sin
(\alpha+\phi)],\\
&  & \sqrt{\frac{C}{\lambda_{1}}}\lambda_{1}r_{2}\cos t[\sin\phi_{1}\cos
(\alpha+\phi)-\cos\phi_{1}\sin(\alpha+\phi)],\\
&  & \sqrt{\frac{C}{\lambda_{2}}}\lambda_{2}r_{1}\sin t[\cos\phi_{2}\cos
(\alpha+\phi)+\sin\phi_{2}\sin(\alpha+\phi)],\\
&  & \sqrt{\frac{C}{\lambda_{2}}}\lambda_{2}r_{1}\sin t[\sin\phi_{2}\cos
(\alpha+\phi)-\cos\phi_{2}\sin(\alpha+\phi)],\\
&  & \frac{C}{2}r_{1}r_{2}\sin(\alpha+\phi))\\
& = & (\sqrt{C\lambda_{1}}r_{2}\cos t\cos(\alpha+\phi_{2}),-\sqrt{C\lambda_{1}%
}r_{2}\cos t\sin(\alpha+\phi_{2}),\\
&  & \sqrt{C\lambda_{2}}r_{1}\sin t\cos(\alpha+\phi_{1}),-\sqrt{C\lambda_{2}%
}r_{1}\sin t\sin(\alpha+\phi_{1}),\\
&  & \frac{C}{2}r_{1}r_{2}\sin(\alpha+\phi))\\
& = & \frac{r_{2}}{2}\sqrt{C\lambda_{1}}\cos t\cos(\alpha+\phi_{2})E_{1}%
+\frac{r_{1}}{2}\sqrt{C\lambda_{2}}\sin t\cos(\alpha+\phi_{1})E_{2}\\
&  & +\frac{r_{2}}{2}\sqrt{C\lambda_{1}}\cos t\sin(\alpha+\phi_{2})E_{3}%
+\frac{r_{1}}{2}\sqrt{C\lambda_{2}}\sin t\sin(\alpha+\phi_{1})E_{4},
\end{array}
\label{eq24}%
\end{equation}
Note $\psi^{2}=\frac{C}{4\lambda_{1}}r_{1}^{2}\sin^{2}t+\frac{C}{4\lambda_{2}%
}r_{2}^{2}\cos^{2}t$, then the associated metric is
\[%
\begin{array}
[c]{c}%
||\frac{\partial\tilde{L}}{\partial t}||^{2}=\frac{C}{4\lambda_{1}}r_{1}%
^{2}\sin^{2}t+\frac{C}{4\lambda_{2}}r_{2}^{2}\cos^{2}t=\psi^{2},\text{
\ }||\frac{\partial\tilde{L}}{\partial s}||^{2}=\lambda_{1}\lambda_{2}\psi
^{2}.
\end{array}
\]
Hence, setting
\[%
\begin{array}
[c]{c}%
A_{1}=\frac{r_{1}}{2\psi}\sqrt{\frac{C}{\lambda_{1}}}\sin t,\text{ \ }%
A_{2}=\frac{r_{2}}{2\psi}\sqrt{\frac{C}{\lambda_{2}}}\cos t,\text{
\ }\Theta_{j}=\alpha+\phi_{j}\text{ }(j=1,2),
\end{array}
\]
so that $A_{1}^{2}+A_{2}^{2}=1$, we choose the orthonormal basis to be%
\begin{equation}%
\begin{array}
[c]{lll}%
e_{1}=\frac{1}{\psi}\frac{\partial\tilde{L}}{\partial t} & = & -A_{1}\cos
\phi_{1}E_{1}+A_{2}\cos\phi_{2}E_{2}+A_{1}\sin\phi_{1}E_{3}-A_{2}\sin\phi
_{2}E_{4},
\end{array}
\label{eq25}%
\end{equation}%
\begin{equation}%
\begin{array}
[c]{lll}%
e_{2}=\frac{1}{\sqrt{\lambda_{1}\lambda_{2}}\psi}\frac{\partial\tilde
{L}}{\partial s} & = & F_{1}E_{1}+F_{2}E_{2}+F_{3}E_{3}+F_{4}E_{4},
\end{array}
\label{eq26}%
\end{equation}
where
\[%
\begin{array}
[c]{c}%
F_{1}=A_{2}\cos\Theta_{2},\text{ }F_{2}=A_{1}\cos\Theta_{1},\text{ }%
F_{3}=A_{2}\sin\Theta_{2},\text{ }F_{4}=A_{1}\sin\Theta_{1},
\end{array}
\]
and%
\begin{equation}%
\begin{array}
[c]{lll}%
\Phi e_{1} & = & -A_{1}\sin\phi_{1}E_{1}+A_{2}\sin\phi_{2}E_{2}-A_{1}\cos
\phi_{1}E_{3}+A_{2}\cos\phi_{2}E_{4},
\end{array}
\label{eq27}%
\end{equation}%
\begin{equation}%
\begin{array}
[c]{lll}%
\Phi e_{2} & = & -F_{3}E_{1}-F_{4}E_{2}+F_{1}E_{3}+F_{2}E_{4}.
\end{array}
\label{eq28}%
\end{equation}
Thus
\[%
\begin{array}
[c]{lll}%
\overline{\nabla}_{e_{1}}e_{1} & = & -\overline{\nabla}_{e_{1}}(\frac{r_{1}%
}{2\psi}\sqrt{\frac{C}{\lambda_{1}}}\sin t\cos\phi_{1}E_{1})+\overline{\nabla
}_{e_{1}}(\frac{r_{2}}{2\psi}\sqrt{\frac{C}{\lambda_{2}}}\cos t\cos\phi
_{2}E_{2})\\
&  & +\overline{\nabla}_{e_{1}}(\frac{r_{1}}{2\psi}\sqrt{\frac{C}{\lambda_{1}%
}}\sin t\sin\phi_{1}E_{3})-\overline{\nabla}_{e_{1}}(\frac{r_{2}}{2\psi}%
\sqrt{\frac{C}{\lambda_{2}}}\cos t\sin\phi_{2}E_{4})\\
& = & -\frac{1}{\psi}\frac{\partial}{\partial t}(\frac{r_{1}}{2\psi}%
\sqrt{\frac{C}{\lambda_{1}}}\sin t\cos\phi_{1}E_{1})-\frac{r_{1}}{2\psi}%
\sqrt{\frac{C}{\lambda_{1}}}\sin t\cos\phi_{1}\overline{\nabla}_{e_{1}}E_{1}\\
&  & +\frac{1}{\psi}\frac{\partial}{\partial t}(\frac{r_{2}}{2\psi}\sqrt
{\frac{C}{\lambda_{2}}}\cos t\cos\phi_{2}E_{2})+\frac{r_{2}}{2\psi}\sqrt
{\frac{C}{\lambda_{2}}}\cos t\cos\phi_{2}\overline{\nabla}_{e_{1}}E_{2}\\
&  & +\frac{1}{\psi}\frac{\partial}{\partial t}(\frac{r_{1}}{2\psi}\sqrt
{\frac{C}{\lambda_{1}}}\sin t\sin\phi_{1}E_{3})+\frac{r_{1}}{2\psi}\sqrt
{\frac{C}{\lambda_{1}}}\sin t\sin\phi_{1}\overline{\nabla}_{e_{1}}E_{3}\\
&  & -\frac{1}{\psi}\frac{\partial}{\partial t}(\frac{r_{2}}{2\psi}\sqrt
{\frac{C}{\lambda_{2}}}\cos t\sin\phi_{2}E_{4})-\frac{r_{2}}{2\psi}\sqrt
{\frac{C}{\lambda_{2}}}\cos t\sin\phi_{2}\overline{\nabla}_{e_{1}}E_{4}\\
& = & -\frac{r_{1}\cos t\cos\phi_{1}}{2\psi^{2}}\sqrt{\frac{C}{\lambda_{1}}%
}[1-\frac{C\sin^{2}t}{4\psi^{2}}(\frac{r_{1}^{2}}{\lambda_{1}}-\frac{r_{2}%
^{2}}{\lambda_{2}})]E_{1}+\frac{Cr_{1}^{2}\sin^{2}t\sin\phi_{1}\cos\phi_{1}%
}{4\lambda_{1}\psi^{2}}\xi\\
&  & -\frac{r_{2}\sin t\cos\phi_{2}}{2\psi^{2}}\sqrt{\frac{C}{\lambda_{2}}%
}[1+\frac{C\cos^{2}t}{4\psi^{2}}(\frac{r_{1}^{2}}{\lambda_{1}}-\frac{r_{2}%
^{2}}{\lambda_{2}})]E_{2}+\frac{Cr_{2}^{2}\cos^{2}t\sin\phi_{2}\cos\phi_{2}%
}{4\lambda_{2}\psi^{2}}\xi\\
&  & +\frac{r_{1}\cos t\sin\phi_{1}}{2\psi^{2}}\sqrt{\frac{C}{\lambda_{1}}%
}[1-\frac{C\sin^{2}t}{4\psi^{2}}(\frac{r_{1}^{2}}{\lambda_{1}}-\frac{r_{2}%
^{2}}{\lambda_{2}})]E_{3}-\frac{Cr_{1}^{2}\sin^{2}t\sin\phi_{1}\cos\phi_{1}%
}{4\lambda_{1}\psi^{2}}\xi\\
&  & +\frac{r_{2}\sin t\sin\phi_{2}}{2\psi^{2}}\sqrt{\frac{C}{\lambda_{2}}%
}[1+\frac{C\cos^{2}t}{4\psi^{2}}(\frac{r_{1}^{2}}{\lambda_{1}}-\frac{r_{2}%
^{2}}{\lambda_{2}})]E_{4}-\frac{Cr_{2}^{2}\cos^{2}t\sin\phi_{2}\cos\phi_{2}%
}{4\lambda_{2}\psi^{2}}\xi\\
& = & -\frac{r_{1}\cos t\cos\phi_{1}}{2\psi^{2}}\sqrt{\frac{C}{\lambda_{1}}%
}[1-\frac{C\sin^{2}t}{4\psi^{2}}(\frac{r_{1}^{2}}{\lambda_{1}}-\frac{r_{2}%
^{2}}{\lambda_{2}})]E_{1}\\
&  & -\frac{r_{2}\sin t\cos\phi_{2}}{2\psi^{2}}\sqrt{\frac{C}{\lambda_{2}}%
}[1+\frac{C\cos^{2}t}{4\psi^{2}}(\frac{r_{1}^{2}}{\lambda_{1}}-\frac{r_{2}%
^{2}}{\lambda_{2}})]E_{2}\\
&  & +\frac{r_{1}\cos t\sin\phi_{1}}{2\psi^{2}}\sqrt{\frac{C}{\lambda_{1}}%
}[1-\frac{C\sin^{2}t}{4\psi^{2}}(\frac{r_{1}^{2}}{\lambda_{1}}-\frac{r_{2}%
^{2}}{\lambda_{2}})]E_{3}\\
&  & +\frac{r_{2}\sin t\sin\phi_{2}}{2\psi^{2}}\sqrt{\frac{C}{\lambda_{2}}%
}[1+\frac{C\cos^{2}t}{4\psi^{2}}(\frac{r_{1}^{2}}{\lambda_{1}}-\frac{r_{2}%
^{2}}{\lambda_{2}})]E_{4},
\end{array}
\]%
\[%
\begin{array}
[c]{lll}%
\overline{\nabla}_{e_{2}}e_{2} & = & \frac{1}{\sqrt{\lambda_{1}\lambda_{2}%
}\psi}\sum_{i=1}^{4}\frac{\partial F_{i}}{\partial s}E_{i}+(F_{1}F_{3}%
+F_{2}F_{4}-F_{3}F_{1}-F_{4}F_{2})\xi\\
& = & \frac{1}{\sqrt{\lambda_{1}\lambda_{2}}\psi}\sum_{i=1}^{4}\frac{\partial
F_{i}}{\partial s}E_{i},
\end{array}
\]
where
\[%
\begin{array}
[c]{c}%
\frac{\partial\psi}{\partial t}=\frac{C}{4\psi}(\frac{r_{1}^{2}}{\lambda_{1}%
}-\frac{r_{2}^{2}}{\lambda_{2}})\sin t\cos t,\text{ \ }\psi_{s}=\frac{C}{4\psi
}r_{1}r_{2}\cos(\alpha+\phi),
\end{array}
\]
and, since $\frac{\partial_{s}A_{j}}{A_{j}}=\frac{r_{j}^{\prime}}{r_{j}}%
-\frac{\psi_{s}}{\psi}$ and $\Theta_{j}^{\prime}=\alpha_{s}+\phi_{j}^{\prime}%
=\alpha_{s}-\frac{\lambda_{j}r_{3-j}}{r_{j}}\sin(\alpha+\phi)$,
\[%
\begin{array}
[c]{c}%
\frac{\partial F_{1}}{\partial s}=(\partial_{s}A_{2})\cos\Theta_{2}%
-A_{2}\Theta_{2}^{\prime}\sin\Theta_{2},\text{ \ }\frac{\partial F_{3}%
}{\partial s}=(\partial_{s}A_{2})\sin\Theta_{2}+A_{2}\Theta_{2}^{\prime}%
\cos\Theta_{2},\\
\frac{\partial F_{2}}{\partial s}=(\partial_{s}A_{1})\cos\Theta_{1}%
-A_{1}\Theta_{1}^{\prime}\sin\Theta_{1},\text{ \ }\frac{\partial F_{4}%
}{\partial s}=(\partial_{s}A_{1})\sin\Theta_{1}+A_{1}\Theta_{1}^{\prime}%
\cos\Theta_{1}.
\end{array}
\]
Hence,%
\begin{equation}%
\begin{array}
[c]{ll}
& \left\langle \overline{\nabla}_{e_{1}}e_{1},\Phi e_{1}\right\rangle \\
= & \frac{Cr_{1}^{2}\sin t\cos t\sin\phi_{1}\cos\phi_{1}}{4\lambda_{1}\psi
^{3}}[1-\frac{C\sin^{2}t}{4\psi^{2}}(\frac{r_{1}^{2}}{\lambda_{1}}-\frac
{r_{2}^{2}}{\lambda_{2}})]\\
& -\frac{Cr_{2}^{2}\sin t\cos t\sin\phi_{2}\cos\phi_{2}}{4\lambda_{2}\psi^{3}%
}[1+\frac{C\cos^{2}t}{4\psi^{2}}(\frac{r_{1}^{2}}{\lambda_{1}}-\frac{r_{2}%
^{2}}{\lambda_{2}})]\\
& -\frac{Cr_{1}^{2}\sin t\cos t\sin\phi_{1}\cos\phi_{1}}{4\lambda_{1}\psi^{3}%
}[1-\frac{C\sin^{2}t}{4\psi^{2}}(\frac{r_{1}^{2}}{\lambda_{1}}-\frac{r_{2}%
^{2}}{\lambda_{2}})]\\
& +\frac{Cr_{2}^{2}\sin t\cos t\sin\phi_{2}\cos\phi_{2}}{4\lambda_{2}\psi^{3}%
}[1+\frac{C\cos^{2}t}{4\psi^{2}}(\frac{r_{1}^{2}}{\lambda_{1}}-\frac{r_{2}%
^{2}}{\lambda_{2}})]\\
= & 0,
\end{array}
\label{eq29}%
\end{equation}%
\begin{equation}%
\begin{array}
[c]{ll}
& \left\langle \overline{\nabla}_{e_{2}}e_{2},\Phi e_{1}\right\rangle \\
= & \frac{1}{\sqrt{\lambda_{1}\lambda_{2}}\psi}\left\{  -A_{1}\left[  \sin
\phi_{1}\frac{\partial F_{1}}{\partial s}+\cos\phi_{1}\frac{\partial F_{3}%
}{\partial s}\right]  +A_{2}\left[  \sin\phi_{2}\frac{\partial F_{2}}{\partial
s}+\cos\phi_{2}\frac{\partial F_{4}}{\partial s}\right]  \right\} \\
= & \frac{1}{\sqrt{\lambda_{1}\lambda_{2}}\psi}\left\{  -A_{1}\left[
(\partial_{s}A_{2})\sin(\alpha+\phi)+A_{2}\Theta_{2}^{\prime}\cos(\alpha
+\phi)\right]  \right. \\
& \left.  +A_{2}\left[  (\partial_{s}A_{1})\sin(\alpha+\phi)+A_{1}\Theta
_{1}^{\prime}\cos(\alpha+\phi)\right]  \right\} \\
= & \frac{A_{1}A_{2}}{\sqrt{\lambda_{1}\lambda_{2}}\psi}\left\{  \left(
\frac{r_{1}^{\prime}}{r_{1}}-\frac{r_{2}^{\prime}}{r_{2}}\right)  \sin
(\alpha+\phi)+\left(  \Theta_{1}^{\prime}-\Theta_{2}^{\prime}\right)
\cos(\alpha+\phi)\right\} \\
= & \frac{A_{1}A_{2}}{\sqrt{\lambda_{1}\lambda_{2}}\psi}\left(  \frac
{\lambda_{1}r_{2}}{r_{1}}-\frac{\lambda_{2}r_{1}}{r_{2}}\right)  \left[
\cos(\alpha+\phi)\sin(\alpha+\phi)-\sin(\alpha+\phi)\cos(\alpha+\phi)\right] \\
= & 0,
\end{array}
\label{eq30}%
\end{equation}
where we used $\sin\phi_{1}\cos\Theta_{2}+\cos\phi_{1}\sin\Theta_{2}%
=\sin(\alpha+\phi)$ and $\cos\phi_{1}\cos\Theta_{2}-\sin\phi_{1}\sin\Theta
_{2}=\cos(\alpha+\phi)$ (and similarly with the indices interchanged),%
\begin{equation}%
\begin{array}
[c]{ll}
& \left\langle \overline{\nabla}_{e_{1}}e_{1},\Phi e_{2}\right\rangle \\
= & \frac{Cr_{1}r_{2}\cos^{2}t}{4\sqrt{\lambda_{1}\lambda_{2}}\psi^{3}}\left[
\cos\phi_{1}\sin(\alpha+\phi_{2})+\sin\phi_{1}\cos(\alpha+\phi_{2})\right]
[1-\frac{C\sin^{2}t}{4\psi^{2}}(\frac{r_{1}^{2}}{\lambda_{1}}-\frac{r_{2}^{2}%
}{\lambda_{2}})]\\
& +\frac{Cr_{1}r_{2}\sin^{2}t}{4\sqrt{\lambda_{1}\lambda_{2}}\psi^{3}}\left[
\cos\phi_{2}\sin(\alpha+\phi_{1})+\sin\phi_{2}\cos(\alpha+\phi_{1})\right]
[1+\frac{C\cos^{2}t}{4\psi^{2}}(\frac{r_{1}^{2}}{\lambda_{1}}-\frac{r_{2}^{2}%
}{\lambda_{2}})]\\
= & \frac{Cr_{1}r_{2}}{4\sqrt{\lambda_{1}\lambda_{2}}\psi^{3}}\sin(\alpha
+\phi),
\end{array}
\label{eq31}%
\end{equation}
and
\begin{equation}%
\begin{array}
[c]{ll}
& \left\langle \overline{\nabla}_{e_{2}}e_{2},\Phi e_{2}\right\rangle \\
= & \frac{1}{\sqrt{\lambda_{1}\lambda_{2}}\psi}\left[  F_{1}\frac{\partial
F_{3}}{\partial s}-F_{3}\frac{\partial F_{1}}{\partial s}+F_{2}\frac{\partial
F_{4}}{\partial s}-F_{4}\frac{\partial F_{2}}{\partial s}\right] \\
= & \frac{1}{\sqrt{\lambda_{1}\lambda_{2}}\psi}\left[  A_{2}^{2}\Theta
_{2}^{\prime}+A_{1}^{2}\Theta_{1}^{\prime}\right] \\
= & \frac{1}{\sqrt{\lambda_{1}\lambda_{2}}\psi}\left[  \alpha_{s}-\frac
{Cr_{1}r_{2}}{4\psi^{2}}\sin(\alpha+\phi)\right] \\
= & \frac{\alpha_{s}}{\sqrt{\lambda_{1}\lambda_{2}}\psi}-\frac{Cr_{1}r_{2}%
}{4\sqrt{\lambda_{1}\lambda_{2}}\psi^{3}}\sin(\alpha+\phi),
\end{array}
\label{eq32}%
\end{equation}
From equations (\ref{eq29}), (\ref{eq30}), (\ref{eq31}) and (\ref{eq32}) the
mean curvature vector $\tilde{H}$ is determined to be%
\begin{equation}%
\begin{array}
[c]{lll}%
\tilde{H} & = & (\overline{\nabla}_{e_{i}}e_{i})^{\perp}\\
& = & \langle\overline{\nabla}_{e_{1}}e_{1},\Phi{e_{1}}\rangle\Phi{e}%
_{1}+\langle\overline{\nabla}_{e_{1}}e_{1},\Phi{e_{2}}\rangle\Phi{e}_{2}\\
&  & +\langle\overline{\nabla}_{e_{2}}e_{2},\Phi{e_{1}}\rangle\Phi{e}%
_{1}+\langle\overline{\nabla}_{e_{2}}e_{2},\Phi{e_{2}}\rangle\Phi{e}_{2}\\
& = & \frac{\alpha_{s}}{\sqrt{\lambda_{1}\lambda_{2}}\psi}\Phi{e_{2}}.
\end{array}
\label{eq33}%
\end{equation}
Using (\ref{eq22}), (\ref{eq27}) and (\ref{eq28}), and writing $B_{1}%
=\frac{r_{1}}{2}\sqrt{\frac{C}{\lambda_{1}}}\cos t$, $B_{2}=\frac{r_{2}}%
{2}\sqrt{\frac{C}{\lambda_{2}}}\sin t$, so that $\tilde{L}=B_{1}\cos\phi
_{1}E_{1}+B_{2}\cos\phi_{2}E_{2}-B_{1}\sin\phi_{1}E_{3}-B_{2}\sin\phi_{2}%
E_{4}+\frac{z}{2}\xi$, we find the normal projection of the position field to
be given by%
\[%
\begin{array}
[c]{lll}%
\left\langle \Phi e_{1},\tilde{L}\right\rangle  & = & -A_{1}B_{1}\sin\phi
_{1}\cos\phi_{1}+A_{2}B_{2}\sin\phi_{2}\cos\phi_{2}\\
&  & +A_{1}B_{1}\cos\phi_{1}\sin\phi_{1}-A_{2}B_{2}\cos\phi_{2}\sin\phi_{2}=0,
\end{array}
\]%
\[%
\begin{array}
[c]{lll}%
\left\langle \Phi e_{2},\tilde{L}\right\rangle  & = & -B_{1}A_{2}\left[
\sin\Theta_{2}\cos\phi_{1}+\cos\Theta_{2}\sin\phi_{1}\right]  -B_{2}%
A_{1}\left[  \sin\Theta_{1}\cos\phi_{2}+\cos\Theta_{1}\sin\phi_{2}\right] \\
& = & -\left(  B_{1}A_{2}+B_{2}A_{1}\right)  \sin(\alpha+\phi)=-\frac
{Cr_{1}r_{2}}{4\sqrt{\lambda_{1}\lambda_{2}}\psi}\sin(\alpha+\phi),
\end{array}
\]
and $\langle X,{\xi}\rangle=\eta(X)=z$. It follows that $\tilde{X}^{\perp}$ is
\begin{equation}%
\begin{array}
[c]{c}%
\tilde{X}^{\perp}=\langle X,\Phi{e_{1}}\rangle\Phi{e_{1}}+\langle
X,\Phi{e_{2}}\rangle\Phi{e_{2}}+\langle X,{\xi}\rangle\xi=-\frac{Cr_{1}r_{2}%
}{4\sqrt{\lambda_{1}\lambda_{2}}\psi}\sin(\alpha+\phi){\Phi}e_{2}+z\xi,
\end{array}
\label{eq34}%
\end{equation}
where we used that $X$ and the Euclidean position field have the same
component in the contact distribution. Hence, (\ref{eq19}), (\ref{eq33}) and
(\ref{eq34}) give
\begin{equation}%
\begin{array}
[c]{lll}%
\widetilde{H}-2\alpha\xi & = & \frac{\alpha_{s}}{\sqrt{\lambda_{1}\lambda_{2}%
}\psi}\Phi{e_{2}}-2\alpha\xi=-\frac{aCr_{1}r_{2}}{4\sqrt{\lambda_{1}%
\lambda_{2}}\psi}\sin(\alpha+\phi)\Phi{e_{2}}+az\xi\\
& = & a\left[  -\frac{Cr_{1}r_{2}}{4\sqrt{\lambda_{1}\lambda_{2}}\psi}%
\sin(\alpha+\phi)\Phi{e_{2}}+z\xi\right] \\
& = & a\tilde{X}^{\perp}.
\end{array}
\label{eq35}%
\end{equation}
Here we used that the last equation of (\ref{eq16}) and (\ref{eq19}) give
\[
\frac{d\alpha}{ds}=-\frac{a}{4}Cr_{1}r_{2}\sin(\alpha+\phi)=-\frac{a}{2}%
\frac{dz}{ds},
\]
that is $-2\alpha=az$ after normalizing the additive constant.

On the other hand, since the frames $e_{1},e_{2}$ are orthonormal, the metric
tensor $g$ is the identity matrix. Thus,
\[%
\begin{array}
[c]{c}%
d\alpha=(\nabla_{e_{1}}\alpha)e_{1}+(\nabla_{e_{2}}\alpha)e_{2}=\frac
{1}{\sqrt{\lambda_{1}\lambda_{2}}\psi}\alpha_{s}e_{2},
\end{array}
\]
and
\[%
\begin{array}
[c]{c}%
\Phi\nabla\alpha=\frac{1}{\sqrt{\lambda_{1}\lambda_{2}}\psi}\alpha_{s}%
\Phi{e_{2}}=\widetilde{H},
\end{array}
\]
so that $\alpha$ is indeed the Legendre angle (\ref{angle}) of $\tilde{L}$.

For case ii), by parameterizing $\sum_{j=1}^{2}\lambda_{j}x_{j}^{2}=C$ with
$x_{1}=\sqrt{\frac{C}{\lambda_{1}}}\sec{t},$ $x_{2}=\sqrt{-\frac{C}%
{\lambda_{1}}}\tan{t}$, we can obtain $\psi^{2}=\frac{Cr_{1}^{2}}{4\lambda
_{1}}\tan^{2}t-\frac{Cr_{2}^{2}}{4\lambda_{2}}\sec^{2}t$. Then, similar to the
proof steps above, we can arrive at the same conclusion. Therefore $\tilde{H}-2\alpha{\xi}=a{X^{\perp}}$, as we wanted to prove.

Finally, as $\phi=\phi_{1}+\phi_{2}$ and $Q(u)=\prod_{j=1}^{2}\left(
\rho_{j}+\lambda_{j}u\right)  $, summing the second equation of (\ref{eq16})
over $j=1,$ $2$ gives
\[%
\begin{array}
[c]{c}%
\frac{d\phi}{ds}=\frac{d\phi_{1}}{ds}+\frac{d\phi_{2}}{ds}=-Q(u)^{1/2}%
\sin(\alpha+\phi)\sum_{j=1}^{2}\frac{\lambda_{j}}{\rho_{j}+\lambda_{j}%
u}=-Q(u)^{1/2}{(\ln{Q(u)})}^{\prime}\sin(\alpha+\phi),
\end{array}
\]
and $\frac{du}{ds}=2Q(u)^{1/2}\cos(\alpha+\phi)$ follows from (\ref{eq15}).
Using above equations and (\ref{eq19}), we have%
\[%
\begin{array}
[c]{ll}
& \frac{d}{ds}Q(u)^{1/2}e^{aCu/8}\sin(\alpha+\phi)\\
= & \frac{1}{2}Q(u)^{1/2}{(\ln{Q(u)})}^{\prime}e^{aCu/8}\sin(\alpha+\phi
)\frac{du}{ds}\\
& +Q(u)^{1/2}e^{aCu/8}[\frac{aC}{8}\sin(\alpha+\phi)\frac{du}{ds}+\cos
(\alpha+\phi)(\frac{d\alpha}{ds}+\frac{d\phi}{ds})]\\
= & Q(u)e^{aCu/8}\sin(\alpha+\phi)\cos(\alpha+\phi)[{(\ln{Q(u)})}^{\prime}%
+\frac{aC}{4}]\\
& -Q(u)e^{aCu/8}\sin(\alpha+\phi)\cos(\alpha+\phi)[{(\ln{Q(u)})}^{\prime}%
+\frac{aC}{4}]\\
= & 0,
\end{array}
\]
so the left-hand side of (\ref{eq20}) is a constant, that is $A\in\mathbb{R}$.
\end{proof}

\section{Classification Results}

In the one-dimensional case, the Lagrangian mean curvature flow is commonly
referred to as the curve shortening flow, which has been extensively studied.
The well-known Abresch-Langer \cite{aj} curve provides a classification under
the assumption that the curve is closed. Furthermore, Halldorsson \cite{h12}
extended this result by removing the closedness condition and also classified
the self-expander case.

Consider the ambient space $(\mathbb{R}^{3}, \Phi, \xi, \eta, g)$ with the
contact 1-form $\eta= \frac{1}{2} dz - \frac{1}{4} (y dx - x dy)$, the Reeb
vector field $\xi= 2 \frac{\partial}{\partial z}$, the associated metric $g =
\frac{1}{4}(dx^{2} + dy^{2}) + \eta\otimes\eta$, and the tensor $\Phi= -dx
\otimes\frac{\partial}{\partial y} + dy \otimes\frac{\partial}{\partial x} +
\frac{x}{2} dx \otimes\frac{\partial}{\partial z} + \frac{y}{2} dy
\otimes\frac{\partial}{\partial z}$.

The calculation in this space is analogous to the $(\mathbb{R}^{5}, \Phi, \xi,
\eta, g)$ space introduced in Section 2. Moreover, we can also prove that it
is a Sasakian space form with $\Phi$-sectional curvature $-3$.

\begin{proof}
[proof of Theorem \ref{A-L}]Under the adapted frame $\{E_{1},E_{2}=\Phi
E_{1},\xi\}$ of $(\mathbb{R}^{3},\Phi,\xi,\eta,g)$, writing $P=y\cos
\alpha+x\sin\alpha$ and $Q=x\cos\alpha-y\sin\alpha$ for the first two
components of $\tilde{\gamma}$, the ODE system gives
\[%
\begin{array}
[c]{lll}%
\tilde{\gamma}^{\prime}(t) & = & \left(  B\cos\alpha(t),-B\sin\alpha
(t),\frac{B}{2}x(t)\right) \\
& = & \frac{B}{2}\left(  \cos\alpha(t)E_{1}+\sin\alpha(t)E_{2}\right)  ,
\end{array}
\]
where the $\xi$-component vanishes, so that $\tilde{\gamma}$ is Legendrian, and
$\Vert\dot{\tilde{\gamma}}(t)\Vert_{g}=\frac{B}{2}$. Hence we get the unit
tangent vector
\[%
\begin{array}
[c]{c}%
e=\frac{\dot{\tilde{\gamma}}(t)}{\Vert\dot{\tilde{\gamma}}(t)\Vert_{g}}%
=\cos\alpha(t)E_{1}+\sin\alpha(t)E_{2},
\end{array}
\]
and the unit normal vector
\[
\hat{N}=\Phi e=-\sin\alpha(t)E_{1}+\cos\alpha(t)E_{2}.
\]
Since $\overline{\nabla}_{e}E_{1}=-\sin\alpha\,\xi$ and $\overline{\nabla}%
_{e}E_{2}=\cos\alpha\,\xi$, it follows that%
\[%
\begin{array}
[c]{lll}%
\overline{\nabla}_{e}e & = & \frac{2}{B}\left[  (\cos\alpha)^{\prime}%
E_{1}+(\sin\alpha)^{\prime}E_{2}\right]  +\cos\alpha\overline{\nabla}_{e}%
E_{1}+\sin\alpha\overline{\nabla}_{e}E_{2}\\
& = & \frac{2x(t)}{B}\left(  -\sin\alpha E_{1}+\cos\alpha E_{2}\right)
=\frac{2x(t)}{B}\hat{N},
\end{array}
\]
so that the mean curvature vector is
\[
H:=\kappa\hat{N}=\overline{\nabla}_{e}^{\perp}e=\frac{2x(t)}{B}\hat{N},\text{
\ \ }\kappa=\frac{2x(t)}{B}.
\]
On the other hand $\nabla\alpha=(e\alpha)e=\frac{2}{B}\alpha^{\prime}%
(t)e=\frac{2x(t)}{B}e$, hence
\[%
\begin{array}
[c]{c}%
\Phi\nabla\alpha=\frac{2x(t)}{B}\hat{N}=H,
\end{array}
\]
i.e. $\alpha(t)$ is the Legendre angle of $\tilde{\gamma}(t)$ in the sense of
(\ref{angle}).

At last, let us check that $\tilde{\gamma}(t)$ satisfies the Legendrian
self-shrinker equation. Under the adapted frame $\tilde{\gamma}=\frac{P}{2}%
E_{1}-\frac{Q}{2}E_{2}+\frac{z}{2}\xi$ with $z=\frac{B}{2}\alpha$, so that
\[%
\begin{array}
[c]{c}%
\left\langle \tilde{\gamma},e\right\rangle =\frac{P\cos\alpha-Q\sin\alpha}%
{2}=\frac{y}{2},\text{ \ \ }\left\langle \tilde{\gamma},\hat{N}\right\rangle
=-\frac{P\sin\alpha+Q\cos\alpha}{2}=-\frac{x}{2},
\end{array}
\]
and, since $X$ and $\tilde{\gamma}$ have the same component in the contact
distribution while $\eta(X)=z$,
\[%
\begin{array}
[c]{lll}%
X^{\perp} & = & \left\langle \tilde{\gamma},\hat{N}\right\rangle \hat{N}%
+\eta(X)\xi=-\frac{x}{2}\hat{N}+\frac{B}{2}\alpha\xi\\
& = & -\frac{B}{4}\left(  H-2\alpha\xi\right)  ,
\end{array}
\]
or
\[
\kappa\hat{N}-2\alpha(t)\xi=-\frac{4}{B}X^{\perp},
\]
as desired.
\end{proof}

In the two-dimensional case, Castro and Lerma \cite{cl} obtained a partial
classification theorem for Lagrangian self-similar solutions. They assumed
that the surface is Hamiltonian stationary, which is equivalent to the angle
function being harmonic. We will generalize their result to the Legendrian
case as follows:

\begin{proof}
[proof of Theorem \ref{c-l}]Let us take the part (4) as an example to show the
theorem, the rest will be similar.

Consider the following model constructed by Castro-Lerma \cite{cl}:
\[%
\begin{array}
[c]{lll}%
\Psi_{v}(e^{is},t) & = & \frac{1}{\sqrt{-2a}}\left(  c_{v}\cos se^{\frac
{it}{s_{v}}},t_{v}\sin se^{is_{v}t}\right)  ,
\end{array}
\]
where we put $s_{v}=\sinh v,c_{v}=\cosh v,t_{v}=\coth v=\frac{c_{v}}{s_{v}}$.

Applying the lifting theorem \ref{P51}, it suffices to solve $\tilde{\Psi}%
_{v}^{\ast}\eta=0$. Writing the two components in polar form $w_{j}%
=R_{j}e^{i\varphi_{j}}$, we have $y_{j}dx_{j}-x_{j}dy_{j}=-R_{j}^{2}%
d\varphi_{j}$, so that
\begin{equation}%
\begin{array}
[c]{lll}%
\tilde{\Psi}_{v}^{\ast}\eta & = & \frac{1}{2}dz+\frac{1}{4}\sum_{j=1}^{2}%
R_{j}^{2}d\varphi_{j}\\
& = & \frac{1}{2}dz+\frac{1}{4}\left(  \frac{c_{v}^{2}}{-2a}\frac{\cos^{2}%
s}{s_{v}}+\frac{t_{v}^{2}s_{v}}{-2a}\sin^{2}s\right)  dt\\
& = & \frac{1}{2}dz-\frac{c_{v}^{2}}{8as_{v}}dt,
\end{array}
\label{eq40}%
\end{equation}
where the last equality uses $t_{v}^{2}s_{v}=\frac{c_{v}^{2}}{s_{v}}$, so that
the two terms combine into the constant $\frac{c_{v}^{2}}{-2as_{v}}(\cos
^{2}s+\sin^{2}s)$. Hence%
\begin{equation}%
\begin{array}
[c]{c}%
dz=\frac{c_{v}^{2}}{4as_{v}}dt.
\end{array}
\label{eq41}%
\end{equation}
Then the Legendrian version is given by $\tilde{\Psi}_{v}:\mathbb{S}^{1}%
\times\mathbb{R}\rightarrow\mathbb{R}^{5},$ $v>0$%
\begin{equation}%
\begin{array}
[c]{c}%
\tilde{\Psi}_{v}(e^{is},t)=\left(  \frac{c_{v}}{\sqrt{-2a}}\cos se^{\frac
{i}{s_{v}}t},\frac{t_{v}}{\sqrt{-2a}}\sin se^{is_{v}t},\frac{c_{v}^{2}%
}{4as_{v}}t+c_{0}\right)  .
\end{array}
\label{eq42}%
\end{equation}
Direct calculation gives its tangent vectors as follows%
\[%
\begin{array}
[c]{c}%
X_{1}=\frac{\partial}{\partial s}\tilde{\Psi}_{v}=\left(  -\frac{c_{v}}%
{\sqrt{-2a}}\sin se^{\frac{i}{s_{v}}t},\frac{t_{v}}{\sqrt{-2a}}\cos
se^{is_{v}t},0\right)  ,
\end{array}
\]%
\[%
\begin{array}
[c]{c}%
X_{2}=\frac{\partial}{\partial t}\tilde{\Psi}_{v}=\left(  \frac{it_{v}}%
{\sqrt{-2a}}\cos se^{\frac{i}{s_{v}}t},\frac{ic_{v}}{\sqrt{-2a}}\sin
se^{is_{v}t},\frac{c_{v}^{2}}{4as_{v}}\right)  ,
\end{array}
\]
where we used $\frac{c_{v}}{s_{v}}=t_{v}$ and $t_{v}s_{v}=c_{v}$; then the
associated metric%
\[%
\begin{array}
[c]{c}%
|X_{1}|^{2}=-\frac{1}{8a}\left(  c_{v}^{2}\sin^{2}s+t_{v}^{2}\cos^{2}s\right)
,
\end{array}
\]%
\[%
\begin{array}
[c]{ccc}%
|X_{2}|^{2} & = & -\frac{1}{8a}\left(  t_{v}^{2}\cos^{2}s+c_{v}^{2}\sin
^{2}s\right)  =|X_{1}|^{2}\triangleq e^{2\varrho(s)},
\end{array}
\]%
\[%
\begin{array}
[c]{c}%
\left\langle X_{1},X_{2}\right\rangle =\frac{1}{4}\mathrm{Re}\left[  \frac
{ic_{v}t_{v}}{2a}\sin s\cos s-\frac{ic_{v}t_{v}}{2a}\sin s\cos s\right]  =0.
\end{array}
\]
Note that $\eta(X_{2})=0$, so that only the transverse part $g^{T}$ of $g$
enters these computations. The Legendre angle is calculated by
\begin{equation}
\tilde{\Psi}_{v}^{\ast}\Omega^{T}=e^{-i\alpha}Vol_{\tilde{\Psi}_{v}%
},\label{eq43}%
\end{equation}
where $\Omega^{T}=\frac{1}{4}dz_{1}\wedge dz_{2}$ is the transverse
holomorphic volume form, normalized so that $|\Omega^{T}|_{g^{T}}=1$. By
straightforward computation, the left side of (\ref{eq43}) is%
\[%
\begin{array}
[c]{lll}%
\tilde{\Psi}_{v}^{\ast}\Omega^{T} & = & \frac{1}{4}\det\nolimits_{\mathbb{C}%
}(X_{1},X_{2})ds\wedge dt\\
& = & \frac{i}{8a}e^{i\left(  s_{v}+\frac{1}{s_{v}}\right)  t}\left(  c_{v}%
^{2}\sin^{2}s+t_{v}^{2}\cos^{2}s\right)  ds\wedge dt\\
& = & -ie^{2\varrho}\cdot e^{i\left(  s_{v}+\frac{1}{s_{v}}\right)  t}ds\wedge
dt,
\end{array}
\]
and the right side of (\ref{eq43}) is%
\[%
\begin{array}
[c]{lll}%
Vol_{\tilde{\Psi}_{v}} & = & \sqrt{\det(g_{v})}ds\wedge dt=|X_{1}%
||X_{2}|ds\wedge dt\\
& = & -\frac{1}{8a}(c_{v}^{2}\sin^{2}s+t_{v}^{2}\cos^{2}s)ds\wedge dt\\
& = & e^{2\varrho}ds\wedge dt.
\end{array}
\]
Since $s_{v}+\frac{1}{s_{v}}=\frac{c_{v}^{2}}{s_{v}}=c_{v}t_{v}$, it follows
\[
\alpha=-c_{v}t_{v}t+C_{1}%
\]
for some constant $C_{1}$.

Similar to the proof of Theorem \ref{Ex} and Theorem \ref{A-L}, we can
calculate $H$ by definition under the adapted frame (\ref{eq2.2}), to see that
$H$ and $\alpha$ satisfy
\[
H^{\#}:=\frac{1}{2}d\eta\lrcorner H=\omega^{T}\lrcorner H=d\alpha,\text{
\ \ that is \ }H=\Phi\nabla\alpha,
\]
where $\omega^{T}=g^{T}(\cdot,\Phi\cdot)=\frac{1}{2}d\eta$ is the transverse
K\"{a}hler form, so that $d\eta\lrcorner H=2d\alpha$. Explicitly, $\nabla
\alpha=-c_{v}t_{v}e^{-2\varrho}X_{2}$ and hence%
\[%
\begin{array}
[c]{c}%
H=-c_{v}t_{v}e^{-2\varrho}\left(  \frac{t_{v}}{\sqrt{-2a}}\cos se^{\frac
{i}{s_{v}}t},\frac{c_{v}}{\sqrt{-2a}}\sin se^{is_{v}t}\right)  .
\end{array}
\]
On the other hand, $\left\langle \tilde{\Psi}_{v},X_{2}\right\rangle =0$ and
$\left\langle \tilde{\Psi}_{v},X_{1}\right\rangle =\frac{t_{v}^{2}-c_{v}^{2}%
}{-8a}\sin s\cos s$, so that, writing $D:=c_{v}^{2}\sin^{2}s+t_{v}^{2}\cos
^{2}s=-8ae^{2\varrho}$,%
\[%
\begin{array}
[c]{lll}%
\tilde{\Psi}_{v}^{\perp} & = & \tilde{\Psi}_{v}-\frac{\left\langle \tilde{\Psi
}_{v},X_{1}\right\rangle }{|X_{1}|^{2}}X_{1}\\
& = & \left(  \frac{t_{v}^{2}}{D}\frac{c_{v}}{\sqrt{-2a}}\cos se^{\frac
{i}{s_{v}}t},\frac{c_{v}^{2}}{D}\frac{t_{v}}{\sqrt{-2a}}\sin se^{is_{v}%
t}\right)  +\frac{z}{2}\xi\\
& = & \frac{1}{8a}H+\left(  \frac{c_{v}t_{v}}{8a}t+C_{2}\right)  \xi
\end{array}
\]
for some constant $C_{2}$, where the identification of the first line with
$\frac{1}{8a}H$ follows from $t_{v}^{2}\frac{c_{v}}{\sqrt{-2a}}=c_{v}%
t_{v}\frac{t_{v}}{\sqrt{-2a}}$ and $c_{v}^{2}\frac{t_{v}}{\sqrt{-2a}}%
=c_{v}t_{v}\frac{c_{v}}{\sqrt{-2a}}$ together with $D=-8ae^{2\varrho}$. Since
$\alpha=-c_{v}t_{v}t+C_{1}$ and $z=\frac{c_{v}^{2}}{4as_{v}}t+c_{0}%
=\frac{c_{v}t_{v}}{4a}t+c_{0}$, the generator (\ref{eqX}) of the contact
dilations satisfies, after normalizing the additive constants,
\[%
\begin{array}
[c]{c}%
\tilde{X}^{\perp}=\tilde{\Psi}_{v}^{\perp}+\frac{z}{2}\xi=\frac{1}{8a}%
H+\frac{c_{v}t_{v}}{4a}t\xi=\frac{1}{8a}\left(  H-2\alpha\xi\right)  ,
\end{array}
\]
and therefore
\begin{equation}
H-2\alpha\xi=8a\tilde{X}^{\perp}~~~~(a<0) \label{eq44}%
\end{equation}

\end{proof}

Smoczyk was the first to prove the non-existence of compact orientable
Lagrangian self-shrinkers with a trivial Maslov form \cite{s4}. As a result,
there are no Lagrangian self-shrinkers with the topology of a sphere. This
result can be extended to the Legendrian case using the Lagrangian projection
and Castro-Lerma's theorem for the Whitney sphere \cite{cl}:

\begin{corollary}
There is no Legendrian self-shrinker topologically a sphere in $(\mathbb{R}%
^{5},\Phi,\xi,\eta,g)$.
\end{corollary}

In 2017, Li and Wang established the following rigidity theorem for the
Clifford torus. Note that for a Legendrian immersion the second fundamental
form takes values in the contact distribution, i.e. $g(A(\cdot,\cdot),\xi)=0$,
so that $\Vert A\Vert_{g}^{2}=\Vert A\Vert_{g^{T}}^{2}$ under the Lagrangian
projection; we shall use the two notations interchangeably. Under some conformal transformations of the ambient space,
their statement can be easily modified as follows:

\begin{proposition}
[Theorem 1.2 of \cite{lw}]\label{lw} Let $f: \Sigma^{2}\to(\mathbb{C}%
^{2},g^{T})$ be a compact orientable Lagrangian self-shrinker
\[%
\begin{array}
[c]{c}%
H=-f^{\perp},
\end{array}
\]
where $g^{T}=\frac{1}{4}\sum_{i=1}^{2}((d{x_{i}})^{2}+(d{y_{i}})^{2})$ is the
transverse part of $g$. If $\|A\|_{g^{T}}^{2}\leq2$, then $\|A\|_{g^{T}}%
^{2}=2$ and $\Sigma^{2}$ is the Clifford torus
\[
f(t,s)=(2e^{it}, 2e^{is}).
\]

\end{proposition}

\begin{proof}
[proof of Theorem \ref{rigidity}]The Lagrangian projection of the mean
curvature vector field of Legendrian immersion coincides with the mean
curvature of the projected (immersed) Lagrangian submanifold. Hence
\[%
\begin{array}
[c]{c}%
\Vert A\Vert_{g^{T}}^{2}\leq2.
\end{array}
\]
Now by the compactness of $\bar{F}$, it is easy to see $f=\Pi\circ F$ is a
compact Lagrangian immersion in $\mathbb{R}^{4}$. It follows from Proposition
\ref{lw} that
\[%
\begin{array}
[c]{c}%
\Vert A\Vert_{g^{T}}^{2}=2
\end{array}
\]
and $f:\mathbb{S}^{1}(1)\times\mathbb{S}^{1}(1)\rightarrow\mathbb{C}^{2}$ is
the Clifford torus $\mathbb{S}^{1}(2)\times\mathbb{S}^{1}(2)$
\[
\Sigma:(t,s)\rightarrow(2e^{it},2e^{is}).
\]
Now we lift this Clifford torus to $(\mathbb{R}^{5},\Phi,\xi,\eta,g)$ by
definition of Legendrian submanifold and compute
\[%
\begin{array}
[c]{c}%
z(x)=\frac{1}{2}\int_{\gamma}\sum_{i}(y_{i}dx_{i}-x_{i}dy_{i})=-2(t+s)
\end{array}
\]
given $z(x_{0})=0$. Hence the lifting result is%
\[
F:(t,s)\rightarrow(2e^{it},2e^{is},-2(t+s)).
\]
By taking $a=-\frac{1}{8}$ in part (2) of Theorem \ref{c-l}, we observe that
it is a Legendrian self-shrinker that satisfies%
\[
H-2\alpha\xi=-X^{\perp}.
\]

On the other hand, it follows from Proposition \ref{P51} that we can lift
$\frac{1}{\sqrt{2}}F(\Sigma)\subset\mathbb{R}^{5}$ to $\bar{F}(\Sigma)\subset
\mathbb{R}^{4}\times\mathbb{S}^{1}:$
\[
\bar{F}:(t,s)\rightarrow(\sqrt{2}e^{it},\sqrt{2}e^{is},e^{-i(t+s)}),
\]
the factor $\frac{1}{\sqrt{2}}$ being the one for which the vertical
coordinate of the Legendrian lift is $z=-(t+s)$, so that $e^{iz}$ is single
valued on $\Sigma$. Thus $\bar{F}(\Sigma)$ is a flat Legendrian torus in
$\mathbb{R}^{4}\times\mathbb{S}^{1}$ and, after normalizing the moduli of the
three complex coordinates, it is the flat minimal generalized Legendrian
Clifford torus
\[
\mathbb{T}_{HL}=\left\{  \tfrac{1}{\sqrt{3}}(e^{i\vartheta_{1}},e^{i\vartheta
_{2}},e^{i\vartheta_{3}}):\vartheta_{1}+\vartheta_{2}+\vartheta_{3}\equiv0\text{
}(\mathrm{mod}\text{ }2\pi)\right\}  \subset\mathbb{S}^{5}%
\]
of \cite{h}, the cone over which is the Harvey-Lawson special Lagrangian cone
$\{|z_{1}|=|z_{2}|=|z_{3}|,\text{ }z_{1}z_{2}z_{3}\in\lbrack0,\infty)\}$ in
$\mathbb{C}^{3}$ (\cite{hl}).
\end{proof}

\end{document}